\RequirePackage{ifpdf}
\ifpdf 
\documentclass[pdftex]{sigma}
\else
\documentclass{sigma}
\fi

\numberwithin{equation}{section}

\begin{document}

\allowdisplaybreaks

\renewcommand{\thefootnote}{$\star$}

\renewcommand{\PaperNumber}{018}

\FirstPageHeading

\ShortArticleName{Inverse Spectral Problems for Tridiagonal $N$ by $N$ Complex Hamiltonians}

\ArticleName{Inverse Spectral Problems\\ for Tridiagonal $\boldsymbol{N}$ by $\boldsymbol{N}$ Complex Hamiltonians\footnote{This paper is a contribution to the Proceedings of the VIIth Workshop ``Quantum Physics with Non-Hermitian Operators''
     (June 29 -- July 11, 2008, Benasque, Spain). The full collection
is available at
\href{http://www.emis.de/journals/SIGMA/PHHQP2008.html}{http://www.emis.de/journals/SIGMA/PHHQP2008.html}}}

\Author{Gusein Sh. GUSEINOV}

\AuthorNameForHeading{G.Sh. Guseinov}

\Address{Department of Mathematics, Atilim University, 06836 Incek, Ankara, Turkey}

\Email{\href{mailto:guseinov@atilim.edu.tr}{guseinov@atilim.edu.tr}}
\URLaddress{\url{http://www.atilim.edu.tr/~guseinov/}}

\ArticleDates{Received November 18, 2008, in f\/inal form February 09,
2009; Published online February 14, 2009}

\Abstract{In this paper, the concept of generalized spectral function is introduced for f\/inite-order tridiagonal symmetric matrices (Jacobi matrices) with complex entries. The structure of the generalized spectral function is described in terms of spectral data consisting of the eigenvalues and normalizing numbers of the matrix. The inverse problems from generalized spectral function as well as from spectral data are investigated. In this way, a procedure for construction of complex tridiagonal matrices having real eigenvalues is obtained.}


\Classification{15A29; 39A10}

\section{Introduction}

Consider the $N\times N$ tridiagonal symmetric matrix\ (Jacobi matrix) with
complex entries%
\begin{equation}
J=\left[
\begin{array}{ccccccc}
b_{0} & a_{0} & 0 & \cdots & 0 & 0 & 0 \\
a_{0} & b_{1} & a_{1} & \cdots & 0 & 0 & 0 \\
0 & a_{1} & b_{2} & \cdots & 0 & 0 & 0 \\
\vdots & \vdots & \vdots & \ddots & \vdots & \vdots & \vdots \\
0 & 0 & 0 & \ldots & b_{N-3} & a_{N-3} & 0 \\
0 & 0 & 0 & \cdots & a_{N-3} & b_{N-2} & a_{N-2} \\
0 & 0 & 0 & \cdots & 0 & a_{N-2} & b_{N-1}%
\end{array}%
\right] ,  \label{1.1}
\end{equation}%
where for each $n,$ $a_{n}$ and $b_{n}$ are arbitrary complex numbers such
that $a_{n}$ is dif\/ferent from zero:%
\begin{equation}
a_{n},b_{n}\in
\mathbb{C}
,\qquad a_{n}\neq 0.  \label{1.2}
\end{equation}

In the real case%
\begin{equation}
a_{n},b_{n}\in
\mathbb{R}
,\qquad a_{n}\neq 0,  \label{1.3}
\end{equation}%
the matrix $J$ is Hermitian (self-adjoint) and in this case many versions of
the inverse spectral problem for $J$ have been investigated in the
literature, see \cite{[1],[2],[3]} and references given therein.

In the complex case (\ref{1.2}), the matrix $J$ is in general non-Hermitian
(non-selfadjoint) and our aim in this paper is to introduce appropriate
spectral data for such a matrix and then consider the inverse spectral
problem consisting in determining of the matrix from its spectral data.

As is known \cite{[4],[5],[6],[7],[8],[9]}, for non-selfadjoint dif\/ferential and
dif\/ference operators a natural spectral characteristic is the so-called
\textit{generalized spectral function} which is a linear continuous
functional on an appropriate linear topological space. In general very
little is known about the structure of generalized spectral functions.

Given the matrix $J$ of the form (\ref{1.1}) with the entries satisfying (\ref{1.2}), consider the eigenvalue problem $Jy=\lambda y$ for a column
vector $y=\{y_{n}\}_{n=0}^{N-1},$ that is equivalent to the second order
linear dif\/ference equation%
\begin{equation}
a_{n-1}y_{n-1}+b_{n}y_{n}+a_{n}y_{n+1}=\lambda y_{n},\qquad n\in
\{0,1,\ldots ,N-1\},\qquad a_{-1}=a_{N-1}=1,  \label{1.4}
\end{equation}%
for $\{y_{n}\}_{n=-1}^{N},$ with the boundary conditions%
\begin{equation}
y_{-1}=y_{N}=0.  \label{1.5}
\end{equation}%
The problem (\ref{1.4}), (\ref{1.5}) is a discrete analogue of the
continuous eigenvalue value problem%
\begin{gather}
\frac{d}{dx}\left[ p(x)\frac{d}{dx}y(x)\right] +q(x)y(x)=\lambda y(x),\qquad x\in \lbrack a,b],  \label{1.6}
\\\
y(a)=y(b)=0,  \label{1.7}
\end{gather}%
where $[a,b]$ is a f\/inite interval.

To the continuous problem
\begin{gather*}
\frac{d}{dx}\left[ p(x)\frac{d}{dx}y(x)\right] +q(x)y(x)=\lambda y(x),\qquad x\in \lbrack 0,\infty ),
\\
y(0)=0,
\end{gather*}
on the semi-inf\/inite interval $[0,\infty )$ there corresponds a Jacobi
matrix of the type (\ref{1.1}) but $J$ being inf\/inite both downwards and to
the right. To the equation in (\ref{1.6}) considered on the whole real axis $(-\infty ,\infty )$ there corresponds a Jacobi matrix which is inf\/inite in
the all four directions: upwards, downwards, to the left, and to the right.

The case of inf\/inite Jacobi matrices was considered earlier in the papers
\cite{[6],[7],[8],[9]} in which the generalized spectral function was
introduced and the inverse problem from the generalized spectral function
was studied. However, in the case of inf\/inite Jacobi matrices the structure
of the generalized spectral function does not allow any explicit description
because of complexity of the structure. Our main achievement in the present
paper is that we describe explicitly the structure of the generalized
spectral function for the f\/inite order Jacobi matrices (\ref{1.1}), (\ref{1.2}).

The paper is organized as follows. In Section~\ref{section2}, the generalized
spectral function is introduced for Jacobi matrices of the form (\ref{1.1})
with the entries satisfying (\ref{1.2}). In Section~\ref{section3}, the inverse
problem from the generalized spectral function is investigated. It turns out
that the matrix~(\ref{1.1}) is not uniquely restored from the generalized
spectral function. There are precisely $2^{N-1}$ distinct Jacobi matrices
possessing the same generalized spectral function. The inverse problem is
solved uniquely from the data consisting of the generalized spectral
function and a sequence $\{\sigma _{1},\sigma _{2},\ldots ,\sigma _{N-1}\}$
of signs $+$ and $-$. Section~\ref{section4} is devoted to some examples. In
Section~\ref{section5}, we describe the structure of the generalized spectral function
and in this way we def\/ine the concept of \textit{spectral data} for matrices
(\ref{1.1}). In Section~\ref{section6}, the inverse problem from the spectral
data is considered. In Section~\ref{section7}, we characterize generalized
spectral functions of real Jacobi matrices among the generalized spectral
functions of complex Jacobi matrices. In Section~\ref{section8}, we describe the
structure of generalized spectral functions and spectral data of real Jacobi
matrices. Finally, in Section~\ref{section9}, we consider inverse problem for
real Jacobi matrices from the spectral data.

Note that considerations of complex (non-Hermitian) Hamiltonians in quantum
mechanics and complex discrete models have recently received a lot of
attention \cite{[10],[11],[12],[13]}. For some recent papers dealing with the
spectral theory of dif\/ference (and dif\/ferential) operators with complex
coef\/f\/icients see \cite{[14],[15],[16],[17],[18]}. Fur further reading on the
spectral theory of the Jacobi dif\/ference equation (three term recurrence
relation) the books \cite{[19],[20],[21],[22],[23]} are excellent sources.

\section{Generalized spectral function}\label{section2}

Given a matrix $J$ of the form (\ref{1.1}) with the entries satisfying (\ref{1.2}). Consider the eigenvalue problem $Jy=\lambda y$ for a column vector $y=\{y_{n}\}_{n=0}^{N-1},$ that is equivalent to the second order linear
dif\/ference equation%
\begin{equation}
a_{n-1}y_{n-1}+b_{n}y_{n}+a_{n}y_{n+1}=\lambda y_{n},\qquad n\in
\{0,1,\ldots ,N-1\},\qquad a_{-1}=a_{N-1}=1,  \label{2.1}
\end{equation}%
for $\{y_{n}\}_{n=-1}^{N},$ with the boundary conditions%
\begin{equation}
y_{-1}=y_{N}=0.  \label{2.2}
\end{equation}%
Denote by $\{P_{n}(\lambda )\}_{n=-1}^{N}$ the solution of equation (\ref{2.1}) satisfying the initial conditions
\begin{equation}
y_{-1}=0,\qquad y_{0}=1.  \label{2.3}
\end{equation}%
Using (\ref{2.3}), we can f\/ind from equation (\ref{2.1}) recurrently the
quantities $P_{n}(\lambda )$ for $n=1,2,\ldots ,N;$ $P_{n}(\lambda )$ is a
polynomial in $\lambda $ of degree $n.$

Thus $\{P_{n}(\lambda )\}_{n=0}^{N}$ is the unique solution of the recursion
relations%
\begin{gather}
b_{0}P_{0}(\lambda )+a_{0}P_{1}(\lambda )  = \lambda P_{0}(\lambda ),  \notag
\\
a_{n-1}P_{n-1}(\lambda )+b_{n}P_{n}(\lambda )+a_{n}P_{n+1}(\lambda )
 = \lambda P_{n}(\lambda ),  \label{2.4}
\\
n\in \{1,2,\ldots ,N-1\},\qquad a_{N-1}=1,\nonumber
\end{gather}%
subject to the initial condition%
\begin{equation}
P_{0}(\lambda )=1.  \label{2.5}
\end{equation}

\begin{lemma}
\label{Lem2.1} The equality%
\begin{equation}
\det \left( J-\lambda I\right) =(-1)^{N}a_{0}a_{1}\cdots
a_{N-2}P_{N}(\lambda )  \label{2.6}
\end{equation}%
holds so that the eigenvalues of the matrix $J$ coincide with the zeros of
the polynomial $P_{N}(\lambda ).$
\end{lemma}

\begin{proof}
To prove (\ref{2.6}) let us set, for each $n\in \{1,2,\ldots ,N\},$%
\begin{equation*}
J_{n}=\left[
\begin{array}{ccccccc}
b_{0} & a_{0} & 0 & \cdots & 0 & 0 & 0 \\
a_{0} & b_{1} & a_{1} & \cdots & 0 & 0 & 0 \\
0 & a_{1} & b_{2} & \cdots & 0 & 0 & 0 \\
\vdots & \vdots & \vdots & \ddots & \vdots & \vdots & \vdots \\
0 & 0 & 0 & \ldots & b_{n-3} & a_{n-3} & 0 \\
0 & 0 & 0 & \cdots & a_{n-3} & b_{n-2} & a_{n-2} \\
0 & 0 & 0 & \cdots & 0 & a_{n-2} & b_{n-1}%
\end{array}%
\right]
\end{equation*}%
and set $\Delta _{n}(\lambda )=\det (J_{n}-\lambda I)$. Note that by $I$ we
denote an identity matrix of needed size. By expanding the determinant $\det
(J_{n+1}-\lambda I)$ by the elements of the last row, it is easy to show that%
\begin{equation*}
\Delta _{n+1}(\lambda )=(b_{n}-\lambda )\Delta _{n}(\lambda
)-a_{n-1}^{2}\Delta _{n-1}(\lambda ),\qquad n=1,2,\ldots ;\qquad
\Delta _{0}(\lambda )=1.
\end{equation*}%
Dividing this equation by the product $a_{0}\cdots a_{n-1},$ we f\/ind that
the sequence%
\begin{equation*}
d_{-1}=0,\text{ }d_{0}=1,\text{ }d_{n}=(-1)^{n}(a_{0}\cdots
a_{n-1})^{-1}\Delta _{n}(\lambda ),\qquad n=1,2,\ldots ,
\end{equation*}%
satisf\/ies (\ref{2.1}), (\ref{2.3}). Then $d_{n}=P_{n}(\lambda ),$ $%
n=0,1,\ldots ,$ and hence we have (\ref{2.6}) because $J_{N}=J$ and $%
a_{N-1}=1.$
\end{proof}

For any nonnegative integer $m$ denote by $
\mathbb{C}
_{m}[\lambda ]$ the ring of all polynomials in $\lambda $ of degree $\leq m$
with complex coef\/f\/icients. A mapping $\Omega :
\mathbb{C}
_{m}[\lambda ]\rightarrow
\mathbb{C}
$ is called a \textit{linear} \textit{functional} if for any $G(\lambda
),H(\lambda )\in
\mathbb{C}
_{m}[\lambda ]$ and $\alpha \in
\mathbb{C}
,$ we have%
\begin{equation*}
\left\langle \Omega ,G(\lambda )+H(\lambda )\right\rangle =\left\langle
\Omega ,G(\lambda )\right\rangle +\left\langle \Omega ,H(\lambda
)\right\rangle \qquad \text{and} \qquad \left\langle \Omega ,\alpha G(\lambda
)\right\rangle =\alpha \left\langle \Omega ,G(\lambda )\right\rangle ,
\end{equation*}%
where $\left\langle \Omega ,G(\lambda )\right\rangle $ denotes the value of $%
\Omega $ on the element (polynomial) $G(\lambda ).$

\begin{theorem}
\label{Th2.2} There exists a unique linear functional $\Omega :%
\mathbb{C}
_{2N}[\lambda ]\rightarrow
\mathbb{C}
$ such that the relations%
\begin{gather}
\left\langle \Omega ,P_{m}(\lambda )P_{n}(\lambda )\right\rangle =\delta
_{mn},\qquad  m,n\in \{0,1,\ldots ,N-1\},  \label{2.7}
\\
\left\langle \Omega ,P_{m}(\lambda )P_{N}(\lambda )\right\rangle =0,\qquad m\in \{0,1,\ldots ,N\},  \label{2.8}
\end{gather}%
hold, where $\delta _{mn}$ is the Kronecker delta.
\end{theorem}

\begin{proof}
First we prove the uniqueness of $\Omega .$ Assume that there exists a
linear functional $\Omega $ possessing the properties (\ref{2.7}) and (\ref{2.8}). The $2N+1$ polynomials%
\begin{equation}
P_{n}(\lambda )\quad (n=0,1,\ldots ,N-1),\qquad P_{m}(\lambda
)P_{N}(\lambda )\quad (m=0,1,\ldots ,N)  \label{2.9}
\end{equation}%
form a basis for the linear space $%
\mathbb{C}
_{2N}[\lambda ]$ because they are linearly independent (their degrees are
distinct) and their number $2N+1=\dim
\mathbb{C}
_{2N}[\lambda ].$ On the other hand, by (\ref{2.7}) and (\ref{2.8}) the
functional $\Omega $ takes on polynomials (\ref{2.9}) completely def\/inite
values:%
\begin{gather}
\left\langle \Omega ,P_{n}(\lambda )\right\rangle =\delta _{0n},\qquad
n\in \{0,1,\ldots ,N-1\},  \label{2.10}
\\
\left\langle \Omega ,P_{m}(\lambda )P_{N}(\lambda )\right\rangle =0,\qquad m\in \{0,1,\ldots ,N\}.  \label{2.11}
\end{gather}%
Therefore $\Omega $ is determined on $%
\mathbb{C}
_{2N}[\lambda ]$ uniquely by the linearity.

To prove existence of $\Omega $ we def\/ine it on the basis polynomials (\ref{2.9}) by (\ref{2.10}), (\ref{2.11}) and then we extend $\Omega $ to over
the whole space $%
\mathbb{C}
_{2N}[\lambda ]$ by linearity. Let us show that the functional $\Omega $
def\/ined in this way satisf\/ies (\ref{2.7}), (\ref{2.8}). Denote%
\begin{equation}
\left\langle \Omega ,P_{m}(\lambda )P_{n}(\lambda )\right\rangle =A_{mn},
\qquad m,n\in \{0,1,\ldots ,N\}.  \label{2.12}
\end{equation}%
Obviously, $A_{mn}=A_{nm}$ for $m,n\in \{0,1,\ldots ,N\}.$ From (\ref{2.10})
and (\ref{2.11}) we have%
\begin{gather}
A_{m0}=A_{0m}=\delta _{m0},\qquad m\in \{0,1,\ldots ,N\},  \label{2.13}
\\
A_{mN}=A_{Nm}=0,\qquad m\in \{0,1,\ldots ,N\}.  \label{2.14}
\end{gather}%
Since $\{P_{n}(\lambda )\}_{0}^{N}$ is the solution of equations (\ref{2.4}), we f\/ind from the f\/irst equation, recalling (\ref{2.5}),
\begin{equation*}
\lambda =b_{0}+a_{0}P_{1}(\lambda ).
\end{equation*}%
Substituting this in the remaining equations of (\ref{2.4}), we obtain%
\begin{gather*}
a_{n-1}P_{n-1}(\lambda )+b_{n}P_{n}(\lambda )+a_{n}P_{n+1}(\lambda
)=b_{0}P_{n}(\lambda )+a_{0}P_{1}(\lambda )P_{n}(\lambda ),
\\
n\in \{1,2,\ldots ,N-1\},\qquad a_{N-1}=1.
\end{gather*}%
Applying the functional $\Omega $ to both sides of the last equation, and
recalling (\ref{2.13}) and (\ref{2.14}), we get%
\begin{equation}
A_{n1}=A_{1n}=\delta _{n1},\qquad n\in \{0,1,\ldots ,N\}.  \label{2.15}
\end{equation}%
Further, since%
\begin{gather*}
a_{m-1}P_{m-1}(\lambda )+b_{m}P_{m}(\lambda )+a_{m}P_{m+1}(\lambda )=\lambda
P_{m}(\lambda ),\qquad m\in \{1,2,\ldots ,N-1\},
\\
a_{n-1}P_{n-1}(\lambda )+b_{n}P_{n}(\lambda )+a_{n}P_{n+1}(\lambda )=\lambda
P_{n}(\lambda ),\qquad n\in \{1,2,\ldots ,N-1\},
\end{gather*}%
we obtain, multiplying the f\/irst of these identities by $P_{n}(\lambda ),$
and multiplying the second by $P_{m}(\lambda ),$ then subtracting the second
result from the f\/irst:%
\begin{gather*}
a_{m-1}P_{m-1}(\lambda )P_{n}(\lambda )+b_{m}P_{m}(\lambda )P_{n}(\lambda
)+a_{m}P_{m+1}(\lambda )P_{n}(\lambda )
\\
\quad{}=a_{n-1}P_{n-1}(\lambda )P_{m}(\lambda )+b_{n}P_{n}(\lambda )P_{m}(\lambda
)+a_{n}P_{n+1}(\lambda )P_{m}(\lambda ),
\qquad\!
m,n\!\in\! \{1,2,\ldots ,N\!-\!1\}.
\end{gather*}%
Applying the functional $\Omega $ to both sides of the last equation, and
recalling (\ref{2.13}), (\ref{2.14}), and~(\ref{2.15}), we obtain for $%
A_{mn} $ the boundary value problem%
\begin{gather}
a_{m-1}A_{m-1,n}+b_{m}A_{mn}+a_{m}A_{m+1,n}=a_{n-1}A_{n-1,m}+b_{n}A_{nm}+a_{n}A_{n+1,m},
\label{2.16}
\\
\qquad{} m,n\in \{1,2,\ldots ,N-1\},\nonumber
\\
A_{n0}=A_{0n}=\delta _{n0},\qquad A_{n1}=A_{1n}=\delta _{n1},\qquad A_{Nn}=A_{nN}=0,  \label{2.17}\\
\qquad {} n\in \{0,1,\ldots ,N\}.\nonumber
\end{gather}%
Using(\ref{2.17}), we can f\/ind from (\ref{2.16}) recurrently all the $A_{mn}$
and the unique solution of problem (\ref{2.16}), (\ref{2.17}) is $%
A_{mn}=\delta _{mn}$ for $m,n\in \{0,1,\ldots ,N-1\}$ and $A_{mN}=0$ for $%
m\in \{0,1,\ldots ,N\}.$
\end{proof}

\begin{definition}
The linear functional $\Omega $ def\/ined in Theorem \ref{Th2.2} we call the
generalized spectral function of the matrix $J$ given in (\ref{1.1}).
\end{definition}

\section{Inverse problem from the generalized spectral function}\label{section3}

The inverse problem is stated as follows:

\begin{enumerate}\itemsep=0pt
\item To see if it is possible to reconstruct the matrix $J,$ given its
generalized spectral function~$\Omega $. If it is possible, to describe the
reconstruction procedure.

\item To f\/ind the necessary and suf\/f\/icient conditions for a given linear
functional $\Omega $ on $%
\mathbb{C}
_{2N}[\lambda ],$ to be the generalized spectral function for some matrix $J$
of the form (\ref{1.1}) with entries belonging to the class (\ref{1.2}).
\end{enumerate}

Since $P_{n}(\lambda )$ is a polynomial of degree $n,$ we can write the
representation%
\begin{equation}
P_{n}(\lambda )=\alpha _{n}\left( \lambda ^{n}+\sum_{k=0}^{n-1}\chi
_{nk}\lambda ^{k}\right) ,\qquad n\in \{0,1,\ldots ,N\}.  \label{3.1}
\end{equation}
Substituting (\ref{3.1}) in (\ref{2.4}), we f\/ind that the coef\/f\/icients $a_{n}$, $b_{n}$ of system (\ref{2.4}) and the quanti\-ties~$\alpha _{n}$, $\chi _{nk}$ of decomposition (\ref{3.1}), are interconnected by the equations%
\begin{gather}
a_{n}=\frac{\alpha _{n}}{\alpha _{n+1}}\quad (0\leq n\leq N-2),\qquad
\alpha _{0}=1,\qquad \alpha _{N}=\alpha _{N-1},  \label{3.2}
\\
b_{n}=\chi _{n,n-1}-\chi _{n+1,n}\quad (0\leq n\leq N-1),\qquad
\chi _{0,-1}=0.  \label{3.3}
\end{gather}

It is easily seen that relations (\ref{2.7}), (\ref{2.8}) are equivalent to
the collection of the relations%
\begin{gather}
\left\langle \Omega ,\lambda ^{m}P_{n}(\lambda )\right\rangle =\frac{\delta
_{mn}}{\alpha _{n}},\qquad m=0,1,\ldots ,n,\qquad n\in \{0,1,\ldots
,N-1\},  \label{3.4}
\\
\left\langle \Omega ,\lambda ^{m}P_{N}(\lambda )\right\rangle =0,\qquad
m=0,1,\ldots ,N.  \label{3.5}
\end{gather}%
In fact, using (\ref{3.1}), we have%
\begin{equation}
\left\langle \Omega ,P_{m}(\lambda )P_{n}(\lambda )\right\rangle =\alpha
_{m}\left\langle \Omega ,\lambda ^{m}P_{n}(\lambda )\right\rangle +\alpha
_{m}\sum_{j=0}^{m-1}\chi _{mj}\left\langle \Omega ,\lambda ^{j}P_{n}(\lambda
)\right\rangle .  \label{3.6}
\end{equation}%
Next, since we have the expansion%
\begin{equation*}
\lambda ^{j}=\sum_{i=0}^{j}c_{i}^{(j)}P_{i}(\lambda ),\qquad j\in
\{0,1,\ldots ,N\},
\end{equation*}%
it follows from (\ref{3.6}) that (\ref{3.4}), (\ref{3.5}) hold if we have (\ref{2.7}), (\ref{2.8}). The converse is also true: if~(\ref{3.4}), (\ref{3.5}) hold, then (\ref{2.7}), (\ref{2.8}) can be obtained from (\ref{3.6}),
in conjunction with (\ref{3.1}).

Let us set%
\begin{equation}
s_{l}= \langle \Omega ,\lambda ^{l} \rangle ,\qquad l\in
\{0,1,\ldots ,2N\},  \label{3.7}
\end{equation}%
that are the ``power moments'' of the functional $\Omega .$

Replacing $P_{n}(\lambda )$ and $P_{N}(\lambda )$ in (\ref{3.4}) and (\ref{3.5}) by their expansions in (\ref{3.1}), we obtain%
\begin{gather}
s_{n+m}+\sum_{k=0}^{n-1}\chi _{nk}s_{k+m}=0,\qquad m=0,1,\ldots ,n-1,\qquad n\in \{1,2,\ldots ,N\},  \label{3.8}
\\
s_{2N}+\sum_{k=0}^{N-1}\chi _{Nk}s_{k+N}=0,  \label{3.9}
\\
s_{2n}+\sum_{k=0}^{n-1}\chi _{nk}s_{k+n}=\frac{1}{\alpha _{n}^{2}},\qquad n\in \{0,1,\ldots ,N-1\}.  \label{3.10}
\end{gather}

Notice that (\ref{3.8}) is the \textit{fundamental equation} of the inverse
problem, in the sense that it enables the problem to be formally solved.
For, if we are given the linear functional $\Omega $ on $%
\mathbb{C}
_{2N}[\lambda ],$ we can f\/ind the quantities $s_{l}$ from (\ref{3.7}) and
then we consider the inhomogeneous system of linear  algebraic equations (%
\ref{3.8}) with unknowns $\chi _{n0},\chi _{n1},\ldots ,\chi _{n,n-1},$ for
every f\/ixed $n\in \{1,2,\ldots ,N\}.$ If this system is uniquely solvable,
and $s_{2n}+\sum\limits_{k=0}^{n-1}\chi _{nk}s_{k+n}\neq 0$ for $n\in \{1,2,\ldots
,N-1\},$ then the entries $a_{n},$ $b_{n}$ of the required matrix $J$ can be
found from (\ref{3.2}) and (\ref{3.3}), respectively, $\alpha _{n}$ being
found from (\ref{3.10}). The next theorem gives the conditions under which
the indicated procedure of solving the inverse problem is rigorously
justif\/ied.

\begin{theorem}
\label{Th3.1}In order for a given linear functional $\Omega $, defined on $%
\mathbb{C}
_{2N}[\lambda ]$, to be the generali\-zed spectral function for some Jacobi
matrix $J$ of the form \eqref{1.1} with entries belonging to the class~\eqref{1.2}, it is necessary and sufficient that the following conditions be
satisfied:
\begin{enumerate}\itemsep=0pt
\item[$(i)$] $\left\langle \Omega ,1\right\rangle =1$  (normalization
condition);

\item[$(ii)$] if, for some polynomial $G(\lambda )$, $\deg G(\lambda )\leq
N-1 $,%
\begin{equation*}
\left\langle \Omega ,G(\lambda )H(\lambda )\right\rangle =0
\end{equation*}%
for all polynomials $H(\lambda )$, $\deg H(\lambda )=\deg G(\lambda )$, then
$G(\lambda )\equiv 0;$

\item[$(iii)$] there exists a polynomial $T(\lambda )$ of degree $N$ such that%
\begin{equation*}
\left\langle \Omega ,G(\lambda )T(\lambda )\right\rangle =0
\end{equation*}%
for all polynomials $G(\lambda )$ with $\deg G(\lambda )\leq N.$
\end{enumerate}
\end{theorem}

\begin{proof}
\textit{Necessity}. We obtain $(i)$ from (\ref{2.7}) with $n=m=0,$ recalling (%
\ref{2.5}). To prove $(ii)$, we write the expansion%
\begin{equation*}
G(\lambda )=\sum_{j=0}^{m}c_{j}^{(m)}P_{j}(\lambda ),\qquad m=\deg
G(\lambda ),
\end{equation*}
and take as $H(\lambda )$ the polynomial%
\begin{equation*}
H(\lambda )=\sum_{j=0}^{m}\overline{c_{j}^{(m)}}P_{j}(\lambda ),
\end{equation*}%
where the bar over a complex number denotes the complex conjugation. Then we
f\/ind from $\left\langle \Omega ,G(\lambda )H(\lambda )\right\rangle =0$
using (\ref{2.7}) that%
\begin{equation*}
\sum_{j=0}^{m}\left\vert c_{j}^{(m)}\right\vert ^{2}=0,
\end{equation*}%
hence $c_{j}^{(m)}=0,$ $j=0,1,\ldots ,m,$ i.e., $G(\lambda )\equiv 0.$ The
statement $(iii)$ of the theorem follows from~(\ref{2.8}) if we take $%
T(\lambda )=P_{N}(\lambda ).$

\textit{Sufficiency.} The proof will be given in several stages.

(a) Given the linear functional $\Omega $, def\/ined on $%
\mathbb{C}
_{2N}[\lambda ]$ and satisfying the conditions of the theorem. Consider
equation (\ref{3.8}) with the unknowns $\chi _{nk},$ $k=0,1,\ldots ,n-1,$ in
which $s_{l}$ are found with the aid of the functional $\Omega $ from
expression (\ref{3.7}). Let us show that this equation has a unique solution
for every f\/ixed $n$ $\in \{1,2,\ldots ,N\}.$ For this, it is suf\/f\/icient to
show that the corresponding homogeneous equation%
\begin{equation}
\sum_{k=0}^{n-1}g_{k}s_{k+m}=0,\qquad m=0,1,\ldots ,n-1,  \label{3.11}
\end{equation}%
has only the zero solution for every $n.$ Assume the contrary. For some $%
n\in \{1,2,\ldots ,N\}$ let equation (\ref{3.11}) have the nonzero solution $%
(g_{k})_{0}^{n-1}.$ Further let $(h_{m})_{0}^{n-1}$ be an arbitrary vector.
We multiply both sides of (\ref{3.11}) by $h_{m}$ and sum over $m$ between $%
0 $ and $n-1;$ we get%
\begin{equation*}
\sum_{m=0}^{n-1}\sum_{k=0}^{n-1}h_{m}g_{k}s_{k+m}=0.
\end{equation*}%
Substituting expression (\ref{3.7}) for $s_{k+m}$ in this equation and
denoting%
\begin{equation*}
G(\lambda )=\sum_{k=0}^{n-1}g_{k}\lambda ^{k},\qquad H(\lambda
)=\sum_{m=0}^{n-1}h_{m}\lambda ^{m},
\end{equation*}%
we obtain%
\begin{equation}
\left\langle \Omega ,G(\lambda )H(\lambda )\right\rangle =0.  \label{3.12}
\end{equation}%
Since $(h_{m})_{0}^{n-1}$ is an arbitrary vector, we f\/ind from (\ref{3.12}),
in the light of condition $(ii)$ of the theorem, that $G(\lambda )\equiv 0,$
and hence $g_{0}=g_{1}=\dots =g_{n-1}=0,$ in spite of our assumption. Thus,
for any $n\in \{1,2,\ldots ,N\},$ equation (\ref{3.8}) has a unique solution.

(b) Let us show that%
\begin{equation}
s_{2n}+\sum_{k=0}^{n-1}\chi _{nk}s_{k+n}\neq 0,\qquad n\in \{1,2,\ldots
,N-1\},  \label{3.13}
\end{equation}%
where $(\chi _{nk})_{k=0}^{n-1}$ is the solution of the fundamental equation
(\ref{3.8}). (For $n=0,$ the left-hand side of (\ref{3.13}) is $%
s_{0}=\left\langle \Omega ,1\right\rangle =1.$) Assume the contrary, i.e.,
for some $n\in \{1,2,\ldots ,N-1\}$%
\begin{equation*}
s_{2n}+\sum_{k=0}^{n-1}\chi _{nk}s_{k+n}=0.
\end{equation*}%
Joining this equation to the fundamental equation (\ref{3.8}), we obtain%
\begin{equation}
s_{n+m}+\sum_{k=0}^{n-1}\chi _{nk}s_{k+m}=0,\qquad m=0,1,\ldots ,n.
\label{3.14}
\end{equation}%
Let $(h_{m})_{0}^{n}$ be an arbitrary vector. Multiplying both sides of (\ref%
{3.14}) by $h_{m}$ and summing over $m$ from $0$ to $n,$ we obtain%
\begin{equation*}
\sum_{m=0}^{n}h_{m}s_{n+m}+\sum_{m=0}^{n}\sum_{k=0}^{n-1}h_{m}\chi
_{nk}s_{k+m}=0.
\end{equation*}%
Replacing $s_{l}$ in this by its expression (\ref{3.7}), we obtain%
\begin{equation*}
\left\langle \Omega ,\left[ \lambda ^{n}+\chi (\lambda )\right] H(\lambda
)\right\rangle =0,
\end{equation*}%
where%
\begin{equation*}
\chi (\lambda )=\sum_{k=0}^{n-1}\chi _{nk}\lambda ^{k},\qquad H(\lambda
)=\sum_{m=0}^{n}h_{m}\lambda ^{m}.
\end{equation*}%
Since $(h_{m})_{0}^{n}$ is an arbitrary vector, we obtain from the last
equation in the light of condition $(ii)$ of the theorem:%
\begin{equation*}
\lambda ^{n}+\chi (\lambda )\equiv 0,
\end{equation*}%
which is impossible. Our assumption is therefore false.

(c) Given the solution $(\chi _{nk})_{k=0}^{n-1}$ of the fundamental equation~(\ref{3.8}), we f\/ind $\alpha _{n}$ from (\ref{3.10}) for $n\in \{0,1,\ldots
,N-1\}$ with $\alpha _{0}=1$ and set $\alpha _{N}=\alpha _{N-1}.$ Then we
f\/ind the polynomials $P_{n}(\lambda )$ from (\ref{3.1}). Let us show that
the relations (\ref{2.7}), (\ref{2.8}) hold. It is enough to show that~(\ref{3.4}),~(\ref{3.5}) hold because (\ref{3.4}), (\ref{3.5}) together are
equivalent to relations (\ref{2.7}), (\ref{2.8}). From (\ref{3.8}) and~(\ref{3.10}) we have (\ref{3.4}) and (\ref{3.5}), the latter except for $m=N.$ So
it remains to show that%
\begin{equation*}
 \langle \Omega ,\lambda ^{N}P_{N}(\lambda ) \rangle =0.
\end{equation*}%
For this purpose we use the condition $(iii)$ of the theorem. By this
condition we have%
\begin{equation*}
T(\lambda )=\sum_{k=0}^{N}t_{k}\lambda ^{k},\qquad t_{N}\neq 0,
\end{equation*}%
and%
\begin{equation*}
0=\left\langle \Omega ,P_{N}(\lambda )T(\lambda )\right\rangle
=\sum_{k=0}^{N}t_{k} \langle \Omega ,\lambda ^{k}P_{N}(\lambda
) \rangle =t_{N} \langle \Omega ,\lambda ^{N}P_{N}(\lambda
) \rangle ,
\end{equation*}%
where we have used (\ref{3.5}) except for $m=N.$ Hence $ \langle \Omega
,\lambda ^{N}P_{N}(\lambda ) \rangle =0.$

(d) Let us show that the polynomials $P_{n}(\lambda ),$ $n=0,1,\ldots ,N,$
constructed in accordance with~(\ref{3.1}) with the aid of the numbers $\chi
_{nk}$ and $\alpha _{n}$ obtained above, satisfy the equations%
\begin{gather}
b_{0}P_{0}(\lambda )+a_{0}P_{1}(\lambda )  = \lambda P_{0}(\lambda ),  \notag
\\
a_{n-1}P_{n-1}(\lambda )+b_{n}P_{n}(\lambda )+a_{n}P_{n+1}(\lambda )
 = \lambda P_{n}(\lambda ),  \label{3.15}
\\
n\in \{1,2,\ldots ,N-1\},\qquad a_{N-1}=1,\nonumber
\end{gather}%
where the coef\/f\/icients $a_{n},$ $b_{n}$ are given by the expressions%
\begin{gather}
a_{n}=\frac{\alpha _{n}}{\alpha _{n+1}}\quad (0\leq n\leq N-2),\qquad \alpha _{0}=1,\qquad \alpha _{N}=\alpha _{N-1},  \label{3.16}
\\
b_{n}=\chi _{n,n-1}-\chi _{n+1,n}\quad (0\leq n\leq N-1),\qquad
\chi _{0,-1}=0.  \label{3.17}
\end{gather}%
We f\/irst verify the f\/irst equation of (\ref{3.15}). From (\ref{3.1}) we have%
\begin{equation*}
P_{0}(\lambda )=1,\qquad P_{1}(\lambda )=\alpha _{1}(\lambda +\chi
_{10}).
\end{equation*}%
Hence the f\/irst equation of (\ref{3.15}) has the form%
\begin{equation*}
b_{0}+a_{0}\alpha _{1}(\lambda +\chi _{10})=\lambda ;
\end{equation*}%
and this is true, since, by (\ref{3.16}) and (\ref{3.17}),%
\begin{equation*}
a_{0}\alpha _{1}=1,\qquad b_{0}=-\chi _{10}.
\end{equation*}%
Let us prove the remaining equations of (\ref{3.15}). Since $\lambda
P_{n}(\lambda )$ is a polynomial of degree $n+1,$ while $P_{k}(\lambda ),$ $%
k=0,1,\ldots ,n+1,$ are linearly independent, we have%
\begin{equation}
\lambda P_{n}(\lambda )=\sum_{k=0}^{n+1}c_{k}^{(n)}P_{k}(\lambda )\qquad n\in \{1,2,\ldots ,N-1\},\qquad c_{N}^{(N-1)}=1,  \label{3.18}
\end{equation}%
where $c_{k}^{(n)},$ $k=0,1,\ldots ,n+1,$ are constants. By (\ref{2.7}), (\ref{2.8}) which we proved in 
(c), we have from (\ref{3.18}):%
\begin{equation}
c_{k}^{(n)}=\left\langle \Omega ,\lambda P_{n}(\lambda )P_{k}(\lambda
)\right\rangle ,\qquad k=0,1,\ldots ,n+1\quad (n\in \{1,2,\ldots
,N-2\}).  \label{3.19}
\end{equation}%
The polynomials $\lambda P_{k}(\lambda ),$ $k=0,1,\ldots ,n-2,$ have degrees
$\leq n-1,$ and hence we f\/ind from~(\ref{3.19}) in the light of (\ref{2.7}),
(\ref{2.8}) that%
\begin{equation*}
c_{k}^{(n)}=0,\qquad k=0,1,\ldots ,n-2\quad (n\in \{1,2,\ldots
,N-1\}).
\end{equation*}%
Consequently, expansion (\ref{3.18}) takes the form%
\begin{equation}
c_{n-1}^{(n)}P_{n-1}(\lambda )+c_{n}^{(n)}P_{n}(\lambda
)+c_{n+1}^{(n)}P_{n+1}(\lambda )=\lambda P_{n}(\lambda ),\qquad n\in
\{1,2,\ldots ,N-1\}.  \label{3.20}
\end{equation}%
It follows from (\ref{3.19}) that $c_{n-1}^{(n)}=c_{n}^{(n-1)}$. Hence,
denoting%
\begin{equation}
c_{n+1}^{(n)}=\widetilde{a}_{n},\qquad c_{n}^{(n)}=\widetilde{b}_{n},
\label{3.21}
\end{equation}%
we have from (\ref{3.20}):%
\begin{equation}
\widetilde{a}_{n-1}P_{n-1}(\lambda )+\widetilde{b}_{n}P_{n}(\lambda )+%
\widetilde{a}_{n}P_{n+1}(\lambda )=\lambda P_{n}(\lambda ),\qquad n\in
\{1,2,\ldots ,N-1\}.  \label{3.22}
\end{equation}%
Replacing $P_{n}(\lambda )$ in (\ref{3.22}) by its expression (\ref{3.1})
and equating coef\/f\/icients of like $\lambda ^{n},$ while recalling (\ref{3.16}%
), (\ref{3.17}), we obtain%
\begin{gather*}
\widetilde{a}_{n}=\frac{\alpha _{n}}{\alpha _{n+1}}=a_{n}\quad (0\leq
n\leq N-2),
\\
\widetilde{b}_{n}=\chi _{n,n-1}-\chi _{n+1,n}=b_{n}\quad (0\leq n\leq
N-1).
\end{gather*}%
Theorem \ref{Th3.1} is completely proved.
\end{proof}

\begin{remark}
It follows from the above solution of the inverse problem that the matrix (%
\ref{1.1}) is not uniquely restored from the generalized spectral function.
This is linked with the fact that the $\alpha _{n}$ are determined from (\ref%
{3.10}) uniquely up to a sign. To ensure that the inverse problem is
uniquely solvable, we have to specify additionally a sequence of signs $+$
and $-.$ Namely, let $\{\sigma _{1},\sigma _{2},\ldots ,\sigma _{N-1}\}$ be
a given f\/inite sequence, where for each $n\in \{1,2,\ldots ,N-1\}$ the $%
\sigma _{n}$ is $+$ or $-.$ We have $2^{N-1}$ such dif\/ferent sequences. Now
to determine $\alpha _{n}$ uniquely from (\ref{3.10}) for $n\in \{1,2,\ldots
,N-1\}$ (remember that we always take $\alpha _{0}=1$) we can choose the
sign $\sigma _{n}$ when extracting the square root. In this way we get
precisely $2^{N-1}$ distinct Jacobi matrices possessing the same generalized
spectral function. For example, the two matrices%
\begin{equation*}
\left[
\begin{array}{cc}
1 & 1 \\
1 & 1%
\end{array}%
\right] ,\qquad \left[
\begin{array}{cc}
1 & -1 \\
-1 & 1%
\end{array}%
\right] ,
\end{equation*}%
as well as the four matrices%
\begin{equation*}
\left[
\begin{array}{ccc}
1 & 1 & 0 \\
1 & 1 & 1 \\
0 & 1 & 1%
\end{array}%
\right] ,\qquad \left[
\begin{array}{ccc}
1 & -1 & 0 \\
-1 & 1 & 1 \\
0 & 1 & 1%
\end{array}%
\right] ,\qquad \left[
\begin{array}{ccc}
1 & 1 & 0 \\
1 & 1 & -1 \\
0 & -1 & 1%
\end{array}%
\right] ,\qquad \left[
\begin{array}{ccc}
1 & -1 & 0 \\
-1 & 1 & -1 \\
0 & -1 & 1%
\end{array}%
\right] ,
\end{equation*}%
have the same generalized spectral function. The inverse problem is solved
uniquely from the data consisting of $\Omega $ and a sequence $\{\sigma
_{1},\sigma _{2},\ldots ,\sigma _{N-1}\}$ of signs $+$ and $-.$
\end{remark}

Using the numbers%
\begin{equation}
s_{l}=\langle \Omega ,\lambda ^{l}\rangle ,\qquad
l=0,1,\ldots ,2N,  \label{3.23}
\end{equation}%
let us introduce the determinants%
\begin{equation}
D_{n}=\left\vert
\begin{array}{cccc}
s_{0} & s_{1} & \cdots & s_{n} \\
s_{1} & s_{2} & \cdots & s_{n+1} \\
\vdots & \vdots & \ddots & \vdots \\
s_{n} & s_{n+1} & \cdots & s_{2n}%
\end{array}%
\right\vert ,\qquad n=0,1,\ldots ,N.  \label{3.24}
\end{equation}%
It turns out that Theorem \ref{Th3.1} is equivalent to the following theorem.

\begin{theorem}
\label{Th3.2}In order for a given linear functional $\Omega $, defined on $%
\mathbb{C}
_{2N}[\lambda ]$, to be the genera\-li\-zed spectral function for some Jacobi
matrix $J$ of the form \eqref{1.1} with entries belonging to the class~\eqref{1.2}, it is necessary and sufficient that
\begin{equation}
D_{0}=1,\qquad D_{n}\neq 0\quad (n=1,2,\ldots ,N-1),\qquad \text{and} \qquad D_{N}=0,  \label{3.25}
\end{equation}%
where $D_{n}$ is defined by \eqref{3.24} and \eqref{3.23}.
\end{theorem}

\begin{proof}
\textit{Necessity}. The condition $D_{0}=1$ follows from $1=\left\langle
\Omega ,1\right\rangle =s_{0}=D_{0}.$ By Theorem \ref{Th3.1}, if for a
polynomial%
\begin{equation}
G(\lambda )=\sum_{k=0}^{n}g_{k}\lambda ^{k}  \label{3.26}
\end{equation}%
with $n\leq N-1$ we have
\begin{equation}
\left\langle \Omega ,G(\lambda )H(\lambda )\right\rangle =0  \label{3.27}
\end{equation}%
for all polynomials
\begin{equation}
H(\lambda )=\sum_{m=0}^{n}h_{m}\lambda ^{m},  \label{3.28}
\end{equation}%
then $G(\lambda )\equiv 0,$ that is, $g_{0}=g_{1}=\dots =g_{n}=0.$

If we substitute (\ref{3.26}) and (\ref{3.28}) in (\ref{3.27}), then we get%
\begin{equation*}
\sum_{m=0}^{n}h_{m}\left( \sum_{k=0}^{n}g_{k}s_{k+m}\right) =0.
\end{equation*}%
Since $h_{0},h_{1},\ldots ,h_{n}$ are arbitrary, the last equation gives%
\begin{equation}
\sum_{k=0}^{n}g_{k}s_{k+m}=0,\qquad m=0,1,\ldots ,n.  \label{3.29}
\end{equation}%
This is a linear homogeneous system of algebraic equations with respect to
$g_{0},g_{1},\ldots ,g_{n}$ and the determinant of this system coincides
with the determinant $D_{n}$. Since this system has only the trivial solution
$g_{0}=g_{1}=\dots =g_{n}=0,$ we have that $D_{n}\neq 0,$ where $n\leq N-1.$

To prove that $D_{N}=0,$ we write equation (\ref{3.8}) for $n=N$ to get%
\begin{equation*}
s_{N+m}+\sum_{k=0}^{N-1}\chi _{Nk}s_{k+m}=0,\qquad m=0,1,\ldots ,N-1.
\end{equation*}%
This equation has the unique solution $\chi _{N0},\chi _{N1},\ldots ,\chi
_{N,N-1}.$ Next, these equalities together with~(\ref{3.9}) can be written
in the form%
\begin{equation*}
\left[
\begin{array}{c}
s_{N} \\
s_{N+1} \\
\vdots \\
s_{2N-1} \\
s_{2N}%
\end{array}%
\right] +\chi _{N0}\left[
\begin{array}{c}
s_{0} \\
s_{1} \\
\vdots \\
s_{N-1} \\
s_{N}%
\end{array}%
\right] +\chi _{N1}\left[
\begin{array}{c}
s_{1} \\
s_{2} \\
\vdots \\
s_{N} \\
s_{N+1}%
\end{array}%
\right] +\dots +\chi _{N,N-1}\left[
\begin{array}{c}
s_{N-1} \\
s_{N} \\
\vdots \\
s_{2N-2} \\
s_{2N-1}%
\end{array}%
\right] =0.
\end{equation*}%
This means that the last column in the determinant $D_{N}$ is a linear
combination of the remaining columns. Therefore $D_{N}=0.$

\textit{Sufficiency}. Given the linear functional $\Omega :%
\mathbb{C}
_{2N}[\lambda ]\rightarrow
\mathbb{C}
$ satisfying the conditions (\ref{3.25}), it is enough to show that then the
conditions of Theorem \ref{Th3.1} are satisf\/ied. We have $\left\langle
\Omega ,1\right\rangle =s_{0}=D_{0}=1.$ Next, let (\ref{3.27}) hold for a
polynomial $G(\lambda )$ of the form (\ref{3.26}) and all polynomials $%
H(\lambda )$ of the form (\ref{3.28}). Then (\ref{3.29}) holds. Since the
determinant of this system is $D_{n}$ and $D_{n}\neq 0$ for $n\leq N-1,$ we
get that $g_{0}=g_{1}=\dots =g_{n}=0,$ that is, $G(\lambda )\equiv 0.$
Finally, we have to show that there is a polynomial $T(\lambda )$ of degree $%
N$ such that%
\begin{equation}
\left\langle \Omega ,G(\lambda )T(\lambda )\right\rangle =0  \label{3.30}
\end{equation}%
for all polynomials $G(\lambda )$ with $\deg G(\lambda )\leq N.$ For this
purpose we consider the homogeneous system%
\begin{equation}
\sum_{k=0}^{N}t_{k}s_{k+m}=0,\qquad m=0,1,\ldots ,N,  \label{3.31}
\end{equation}%
with the unknowns $t_{0},t_{1},\ldots ,t_{N}.$ The determinant of this
system is $D_{N}.$ Since by condition $D_{N}=0,$ this system has a
nontrivial solution $t_{0},t_{1},\ldots ,t_{N}.$ We have $t_{N}\neq 0.$
Indeed, if $t_{N}=0,$ then we get from (\ref{3.31})%
\begin{equation}
\sum_{k=0}^{N-1}t_{k}s_{k+m}=0,\qquad m=0,1,\ldots ,N-1.  \label{3.32}
\end{equation}%
The determinant of this system is $D_{N-1}$ and by condition $D_{N-1}\neq 0.$
Then $t_{0}=t_{1}=\dots =t_{N-1}=0$ and we get that the solution $%
t_{0},t_{1},\ldots ,t_{N}$ of system (\ref{3.31}) is trivial, which is a~contradiction. Taking the nontrivial solution $t_{0},t_{1},\ldots ,t_{N}$ of
system (\ref{3.31}) we construct the polynomial%
\begin{equation*}
T(\lambda )=\sum_{k=0}^{N}t_{k}\lambda ^{k}.
\end{equation*}%
of degree $N.$ Then substituting $s_{k+m}= \langle \Omega ,\lambda
^{k+m} \rangle $ in (\ref{3.31}) gives%
\begin{equation*}
\left\langle \Omega ,\lambda ^{m}T(\lambda )\right\rangle =0,\qquad
m=0,1,\ldots ,N.
\end{equation*}%
Hence (\ref{3.30}) holds for all polynomials $G(\lambda )$ with $\deg
G(\lambda )\leq N.$
\end{proof}

Note that the determinant of system (\ref{3.8}) coincides with $D_{n-1}.$
Denote by $D_{m}^{(k)}$ $(k=0,1,\ldots ,m)$ the determinant that is obtained
from the determinant $D_{m}$ by replacing in $D_{m}$ the $(k+1)$th column by
the column with the components $s_{m+1},s_{m+2},\ldots ,s_{2m+1}.$ Then,
solving system (\ref{3.8}) by making use of the Cramer's rule, we f\/ind
\begin{equation}
\chi _{nk}=-\frac{D_{n-1}^{(k)}}{D_{n-1}},\qquad k=0,1,\ldots ,n-1.
\label{3.33}
\end{equation}%
Next, substituting the expression (\ref{3.33}) of $\chi _{nk}$ into the
left-hand side of (\ref{3.10}), we get%
\begin{equation}
\alpha _{n}^{-2}=D_{n}D_{n-1}^{-1}.  \label{3.34}
\end{equation}%
Now if we set $D_{m}^{(m)}=\Delta _{m},$ then we get from (\ref{3.16}), (\ref{3.17}), by virtue of (\ref{3.33}), (\ref{3.34}),
\begin{gather}
a_{n}=\pm \left( D_{n-1}D_{n+1}\right) ^{1/2}D_{n}^{-1},\qquad n\in
\{0,1,\ldots ,N-2\},\qquad D_{-1}=1,  \label{3.35}
\\
b_{n}=\Delta _{n}D_{n}^{-1}-\Delta _{n-1}D_{n-1}^{-1},\qquad n\in
\{0,1,\ldots ,N-1\},\qquad \Delta _{-1}=0,\qquad \Delta _{0}=s_{1}.
\label{3.36}
\end{gather}

Thus, if the conditions of Theorem \ref{Th3.2} or, equivalently, the
conditions of Theorem \ref{Th3.1} are sa\-tisf\/ied, then the entries $a_{n},$ $%
b_{n}$ of the matrix $J$ for which $\Omega $ is the generalized spectral
function, are recovered by the formulas (\ref{3.35}), (\ref{3.36}), where $%
D_{n}$ is def\/ined by (\ref{3.24}) and (\ref{3.23}), and $\Delta _{n}$ is the
determinant obtained from the determinant $D_{n}$ by replacing in $D_{n}$
the last column by the column with the components $s_{n+1},s_{n+2},\ldots
,s_{2n+1}.$

\section{Examples}\label{section4}

In this section we consider some simple examples to illustrate the solving
of the inverse problem given above in Section~\ref{section3}.

\begin{example}
The functional%
\begin{equation*}
\left\langle \Omega ,G(\lambda )\right\rangle =\int_{0}^{1}G(\lambda
)d\lambda
\end{equation*}%
satisf\/ies the conditions $(i)$ and $(ii)$ of Theorem \ref{Th3.1}, but it does
not satisfy the condition $(iii)$ of this theorem.
\end{example}

In fact, obviously, $\langle \Omega ,1 \rangle =1.$ Next, let for
a polynomial%
\begin{equation}
G(\lambda )=\sum_{k=0}^{N-1}g_{k}\lambda ^{k}  \label{4.1}
\end{equation}%
we have%
\begin{equation*}
\left\langle \Omega ,G(\lambda )H(\lambda )\right\rangle
=\int_{0}^{1}G(\lambda )H(\lambda )d\lambda =0
\end{equation*}%
for all polynomials%
\begin{equation}
H(\lambda )=\sum_{k=0}^{N-1}h_{k}\lambda ^{k},\qquad \deg H(\lambda
)=\deg G(\lambda ).  \label{4.2}
\end{equation}%
Taking, in particular,%
\begin{equation}
H(\lambda )=\sum_{k=0}^{N-1}\overline{g}_{k}\lambda ^{k},  \label{4.3}
\end{equation}%
where the bar over a complex number denotes the complex conjugation, we get%
\begin{equation*}
\int_{0}^{1}\left\vert G(\lambda )\right\vert ^{2}d\lambda =0
\end{equation*}%
and hence $G(\lambda )\equiv 0.$

The same reasoning shows that there is no nonidentically zero polynomial $%
T(\lambda )$ such that $\left\langle \Omega ,G(\lambda )T(\lambda
)\right\rangle =0$ for all polynomials $G(\lambda )$ with $\deg G(\lambda
)\leq \deg H(\lambda ).$

\begin{example}
The functional%
\begin{equation*}
\left\langle \Omega ,G(\lambda )\right\rangle =\sum_{k=1}^{N}c_{k}G(\lambda
_{k}),
\end{equation*}%
where $\lambda _{1},\ldots ,\lambda _{N}$ are distinct real numbers, $%
c_{1},\ldots ,c_{N}$ are complex numbers such that%
\begin{equation*}
\sum_{k=1}^{N}c_{k}=1\qquad \text{and} \qquad \text{Re}\; c_{k}>0\quad (k=1,\ldots
,N),
\end{equation*}%
satisf\/ies the conditions of Theorem \ref{Th3.1}.
\end{example}

In fact, obviously, $\left\langle \Omega ,1\right\rangle =1.$ Next, assume
that for a polynomial $G(\lambda )$ of the form (\ref{4.1}) we have%
\begin{equation*}
\left\langle \Omega ,G(\lambda )H(\lambda )\right\rangle
=\sum_{k=1}^{N}c_{k}G(\lambda _{k})H(\lambda _{k})=0
\end{equation*}%
for all polynomials $H(\lambda )$ of the form (\ref{4.2}). If we take, in
particular, $H(\lambda )$ of the form (\ref{4.3}), then we get%
\begin{equation*}
\sum_{k=1}^{N}c_{k}\left\vert G(\lambda _{k})\right\vert ^{2}=0.
\end{equation*}%
Hence, taking the real part and using the condition Re $c_{k}>0$ $%
(k=1,\ldots ,N),$ we get $G(\lambda _{k})=0,$ $k=1,\ldots ,N.$ Therefore $%
G(\lambda )\equiv 0$ because $\lambda _{1},\ldots ,\lambda _{N}$ are
distinct and $G(\lambda )$ is a polynomial with $\deg G(\lambda )\leq N-1.$

Further, for the polynomial%
\begin{equation*}
T(\lambda )=(\lambda -\lambda _{1})\cdots (\lambda -\lambda _{N})
\end{equation*}%
we have $\left\langle \Omega ,G(\lambda )T(\lambda )\right\rangle =0$ for
all polynomials $G(\lambda )$ so that the condition $(iii)$ of Theorem \ref{Th3.1} is also satisf\/ied. Thus the functional $\Omega $ satisf\/ies all the
conditions of Theorem \ref{Th3.1}.

Consider the case $N=2$ and take the functional $\Omega $ def\/ined by the
formula%
\begin{equation*}
\left\langle \Omega ,G(\lambda )\right\rangle =cG(0)+(1-c)G(1),
\end{equation*}%
where $\ c$ is any complex number such that $c\neq 0$ and $c\neq 1.$ Let us
solve the inverse problem for this functional by using formulas (\ref{3.35})
and (\ref{3.36}). We have%
\begin{gather*}
s_{0}= \langle \Omega ,1 \rangle =1,\qquad s_{l}= \langle
\Omega ,\lambda ^{l} \rangle =1-c\qquad \text{for all} \quad l=1,2,\ldots ,
\\
D_{-1}=1,\qquad D_{0}=s_{0}=1,
\\
D_{1}=\left\vert
\begin{array}{cc}
s_{0} & s_{1} \\
s_{1} & s_{2}%
\end{array}%
\right\vert =\left\vert
\begin{array}{cc}
1 & 1-c \\
1-c & 1-c%
\end{array}%
\right\vert =c(1-c),
\\
D_{2}=\left\vert
\begin{array}{ccc}
s_{0} & s_{1} & s_{2} \\
s_{1} & s_{2} & s_{3} \\
s_{2} & s_{3} & s_{4}%
\end{array}%
\right\vert =\left\vert
\begin{array}{ccc}
1 & 1-c & 1-c \\
1-c & 1-c & 1-c \\
1-c & 1-c & 1-c%
\end{array}%
\right\vert =0,
\\
\Delta _{-1}=0,\qquad \Delta _{0}=s_{1}=1-c,
\\
\Delta _{1}=D_{1}^{(1)}=\left\vert
\begin{array}{cc}
s_{0} & s_{2} \\
s_{1} & s_{3}%
\end{array}%
\right\vert =\left\vert
\begin{array}{cc}
1 & 1-c \\
1-c & 1-c%
\end{array}%
\right\vert =c(1-c).
\end{gather*}%
Therefore the functional $\Omega $ satisf\/ies all the conditions of Theorem %
\ref{Th3.2}. According to formu\-las~(\ref{3.35}) and (\ref{3.36}), we f\/ind%
\begin{gather*}
a_{0}=\pm (D_{-1}D_{1})^{1/2}D_{0}^{-1}=\pm \sqrt{D_{1}}=\pm \sqrt{c(1-c)},
\\
b_{0}=\Delta _{0}D_{0}^{-1}-\Delta _{-1}D_{-1}^{-1}=1-c,
\\
b_{1}=\Delta _{1}D_{1}^{-1}-\Delta _{0}D_{0}^{-1}=1-(1-c)=c.
\end{gather*}%
Therefore there are two matrices $J_{\pm }$ for which $\Omega $ is the
spectral function:%
\begin{equation*}
J_{\pm }=\left[
\begin{array}{cc}
b_{0} & a_{0} \\
a_{0} & b_{1}%
\end{array}%
\right] =\left[
\begin{array}{cc}
1-c & \pm \sqrt{c(1-c)} \\
\pm \sqrt{c(1-c)} & c%
\end{array}%
\right] .
\end{equation*}%
The characteristic polynomials of the matrices $J_{\pm }$ have the form%
\begin{equation*}
\det (J_{\pm }-\lambda I)=\lambda (\lambda -1).
\end{equation*}

\begin{example}
Let $N=2.$ Consider the functional $\Omega $ def\/ined by the formula%
\begin{equation*}
\left\langle \Omega ,G(\lambda )\right\rangle =G(\lambda _{0})+cG^{\prime
}(\lambda _{0}),
\end{equation*}%
where $\lambda _{0}$ and $c$ are arbitrary complex numbers such that $c\neq
0.$ This functional satisf\/ies all the conditions of Theorem~\ref{Th3.1}. As
the polynomial $T(\lambda )$ presented in the condition $(iii)$ of Theorem~\ref{Th3.1}, we can take $T(\lambda )=(\lambda -\lambda _{0})^{2}$.
\end{example}

We have%
\begin{gather*}
s_{0}=\left\langle \Omega ,1\right\rangle =1,\qquad s_{l}= \langle
\Omega ,\lambda ^{l} \rangle =\lambda _{0}^{l}+cl\lambda _{0}^{l-1}%
\qquad \text{for} \quad l=1,2,\ldots ,
\\
D_{-1}=1,\qquad D_{0}=s_{0}=1,
\\
D_{1}=\left\vert
\begin{array}{cc}
s_{0} & s_{1} \\
s_{1} & s_{2}%
\end{array}%
\right\vert =\left\vert
\begin{array}{cc}
1 & \lambda _{0}+c \\
\lambda _{0}+c & \lambda _{0}^{2}+2c\lambda _{0}%
\end{array}%
\right\vert =-c^{2},
\\
D_{2}=\left\vert
\begin{array}{ccc}
s_{0} & s_{1} & s_{2} \\
s_{1} & s_{2} & s_{3} \\
s_{2} & s_{3} & s_{4}%
\end{array}%
\right\vert =\left\vert
\begin{array}{ccc}
1 & \lambda _{0}+c & \lambda _{0}^{2}+2c\lambda _{0} \\
\lambda _{0}+c & \lambda _{0}^{2}+2c\lambda _{0} & \lambda
_{0}^{3}+3c\lambda _{0}^{2} \\
\lambda _{0}^{2}+2c\lambda _{0} & \lambda _{0}^{3}+3c\lambda _{0}^{2} &
\lambda _{0}^{4}+4c\lambda _{0}^{3}%
\end{array}%
\right\vert =0,
\\
\Delta _{-1}=0,\qquad \Delta _{0}=s_{1}=\lambda _{0}+c,
\\
\Delta _{1}=D_{1}^{(1)}=\left\vert
\begin{array}{cc}
s_{0} & s_{2} \\
s_{1} & s_{3}%
\end{array}%
\right\vert =\left\vert
\begin{array}{cc}
1 & \lambda _{0}^{2}+2c\lambda _{0} \\
\lambda _{0}+c & \lambda _{0}^{3}+3c\lambda _{0}^{2}%
\end{array}%
\right\vert =-2c^{2}\lambda _{0}.
\end{gather*}
Therefore the functional $\Omega $ satisf\/ies all the conditions of Theorem %
\ref{Th3.2}. According to formu\-las~(\ref{3.35}) and (\ref{3.36}), we f\/ind%
\begin{gather*}
a_{0}=\pm (D_{-1}D_{1})^{1/2}D_{0}^{-1}=\pm \sqrt{D_{1}}=\pm \sqrt{-c^{2}}%
=\pm ic,
\\
b_{0}=\Delta _{0}D_{0}^{-1}-\Delta _{-1}D_{-1}^{-1}=\lambda _{0}+c,
\\
b_{1}=\Delta _{1}D_{1}^{-1}-\Delta _{0}D_{0}^{-1}=\frac{-2c^{2}\lambda _{0}}{%
-c^{2}}-(\lambda _{0}+c)=\lambda _{0}-c.
\end{gather*}%
Therefore the two matrices $J_{\pm }$ for which $\Omega $ is the spectral
function have the form%
\begin{equation*}
J_{\pm }=\left[
\begin{array}{cc}
b_{0} & a_{0} \\
a_{0} & b_{1}%
\end{array}%
\right] =\left[
\begin{array}{cc}
\lambda _{0}+c & \pm ic \\
\pm ic & \lambda _{0}-c%
\end{array}%
\right] .
\end{equation*}%
The characteristic polynomials of the matrices $J_{\pm }$ have the form%
\begin{equation*}
\det (J_{\pm }-\lambda I)=(\lambda -\lambda _{0})^{2}.
\end{equation*}

Note that if $N=3,$ then the functional%
\begin{equation*}
\left\langle \Omega ,G(\lambda )\right\rangle =G(\lambda
_{0})+c_{1}G^{\prime }(\lambda _{0})+c_{2}G^{\prime \prime }(\lambda _{0}),
\end{equation*}%
where $\lambda _{0},$ $c_{1},$ $c_{2}$ are complex numbers, satisf\/ies the
conditions of Theorem \ref{Th3.1} (or Theorem \ref{Th3.2}) if and only if%
\begin{equation*}
c_{2}\neq 0,\qquad 2c_{2}-c_{1}^{2}\neq 0.
\end{equation*}

\section[Structure of the generalized spectral function and spectral data]{Structure of the generalized spectral function\\ and spectral data}\label{section5}

Let $J$ be a Jacobi matrix of the form (\ref{1.1}) with the entries
satisfying (\ref{1.2}). Next, let $\Omega $ be the generalized spectral
function of $J,$ def\/ined above in Section~\ref{section2}. The following theorem
describes the structure of $\Omega .$

\begin{theorem}
\label{Th5.1} Let $\lambda _{1},\ldots ,\lambda _{p}$ be all the distinct
eigenvalues of the matrix $J$ and $m_{1},\ldots ,m_{p}$ be their
multiplicities, respectively, as roots of the characteristic polynomial \eqref{2.6}. There exist complex numbers $\beta _{kj}$ $(j=1,\ldots ,m_{k},$ $k=1,\ldots ,p)$ uniquely determined by the matrix $J$ such that for any
polynomial $G(\lambda )\in
\mathbb{C}
_{2N}[\lambda ]$ the formula%
\begin{equation}
\left\langle \Omega ,G(\lambda )\right\rangle
=\sum_{k=1}^{p}\sum_{j=1}^{m_{k}}\frac{\beta _{kj}}{(j-1)!}G^{(j-1)}(\lambda
_{k}),  \label{5.1}
\end{equation}%
holds, where $G^{(n)}(\lambda )$ denotes the $n$th order derivative of $G(\lambda )$ with respect to $\lambda .$
\end{theorem}

\begin{proof}
Let $J$ be a matrix of the form (\ref{1.1}), (\ref{1.2}). Consider the
second order linear dif\/ference equation%
\begin{equation}
a_{n-1}y_{n-1}+b_{n}y_{n}+a_{n}y_{n+1}=\lambda y_{n},\qquad n\in
\{0,1,\ldots ,N-1\},\qquad a_{-1}=a_{N-1}=1,  \label{5.2}
\end{equation}%
where $\{y_{n}\}_{n=-1}^{N}$ is a desired solution. Denote by $%
\{P_{n}(\lambda )\}_{n=-1}^{N}$ and $\{Q_{n}(\lambda )\}_{n=-1}^{N}$ the
solutions of equation (\ref{5.2}) satisfying the initial conditions%
\begin{gather}
P_{-1}(\lambda )=0,\qquad P_{0}(\lambda )=1;  \label{5.3}
\\
Q_{-1}(\lambda )=-1,\qquad Q_{0}(\lambda )=0.  \label{5.4}
\end{gather}%
For each $n\geq 0,$ $P_{n}(\lambda )$ is a polynomial of degree $n$ and is
called a polynomial of f\/irst kind (note that $P_{n}(\lambda )$ is the same
polynomial as in Section~\ref{section2}) and $Q_{n}(\lambda )$ is a polynomial
of degree $n-1$ and is known as a polynomial of second kind.

Let us set
\begin{equation}
M(\lambda )=-\frac{Q_{N}(\lambda )}{P_{N}(\lambda )}.  \label{5.5}
\end{equation}%
Then it is straightforward to verify that the entries $R_{nm}(\lambda )$ of
the matrix $R(\lambda )=(J-\lambda I)^{-1}$ (resolvent of $J$) are of the
form%
\begin{equation}
R_{nm}(\lambda )=\left\{
\begin{array}{c}
P_{n}(\lambda )[Q_{m}(\lambda )+M(\lambda )P_{m}(\lambda )],\qquad
0\leq n\leq m\leq N-1, \vspace{1mm}\\
P_{m}(\lambda )[Q_{n}(\lambda )+M(\lambda )P_{n}(\lambda )],\qquad
0\leq m\leq n\leq N-1.%
\end{array}%
\right.  \label{5.6}
\end{equation}%
Let $f$ be an arbitrary element (column vector) of $%
\mathbb{C}
^{N},$ with the components $f_{0},f_{1},\ldots ,f_{N-1}.$ Since%
\begin{equation*}
R(\lambda )f=-\frac{f}{\lambda }+O\left( \frac{1}{\lambda ^{2}}\right) ,
\end{equation*}%
as $\left\vert \lambda \right\vert \rightarrow \infty ,$ we have for each $%
n\in \{0,1,\ldots ,N-1\},$%
\begin{equation}
f_{n}=-\frac{1}{2\pi i}\int\nolimits_{\Gamma _{r}}\left\{
\sum_{m=0}^{N-1}R_{nm}(\lambda )f_{m}\right\} d\lambda
+\int\nolimits_{\Gamma _{r}}O\left( \frac{1}{\lambda ^{2}}\right) d\lambda
,  \label{5.7}
\end{equation}%
where $\ r$ is a suf\/f\/iciently large positive number, $\Gamma _{r}$ is the
circle in the $\lambda $-plane of radius $r$ centered at the origin.

Denote by $\lambda _{1},\ldots ,\lambda _{p}$ all the distinct roots of the
polynomial $P_{N}(\lambda )$ (which coincides by (\ref{2.6}) with the
characteristic polynomial of the matrix $J$ up to a constant factor) and by $%
m_{1},\ldots ,m_{p}$ their multiplicities, respectively:%
\begin{equation}
P_{N}(\lambda )=c(\lambda -\lambda _{1})^{m_{1}}\cdots (\lambda -\lambda
_{p})^{m_{p}},  \label{5.8}
\end{equation}%
where $c$ is a constant. We have $1\leq p\leq N$ and $m_{1}+\dots +m_{p}=N.$
By (\ref{5.8}), we can rewrite the rational function $Q_{N}(\lambda
)/P_{N}(\lambda )$ as the sum of partial fractions:%
\begin{equation}
\frac{Q_{N}(\lambda )}{P_{N}(\lambda )}=\sum_{k=1}^{p}\sum_{j=1}^{m_{k}}%
\frac{\beta _{kj}}{(\lambda -\lambda _{k})^{j}},  \label{5.9}
\end{equation}%
where $\beta _{kj}$ are some uniquely determined complex numbers depending
on the matrix $J.$ Substituting (\ref{5.6}) in (\ref{5.7}) and taking into
account (\ref{5.5}), (\ref{5.9}) we get, applying the residue theorem and
passing then to the limit as $r\rightarrow \infty ,$%
\begin{equation}
f_{n}=\sum_{k=1}^{p}\sum_{j=1}^{m_{k}}\frac{\beta _{kj}}{(j-1)!}\left\{
\frac{d^{j-1}}{d\lambda ^{j-1}}\left[ F(\lambda )P_{n}(\lambda )\right]
\right\} _{\lambda =\lambda _{k}},\qquad n\in \{0,1,\ldots ,N-1\},
\label{5.10}
\end{equation}%
where%
\begin{equation}
F(\lambda )=\sum_{m=0}^{N-1}f_{m}P_{m}(\lambda ).  \label{5.11}
\end{equation}%
Now def\/ine on $%
\mathbb{C}
_{2N}[\lambda ]$ the functional $\Omega $ by the formula%
\begin{equation}
\left\langle \Omega ,G(\lambda )\right\rangle
=\sum_{k=1}^{p}\sum_{j=1}^{m_{k}}\frac{\beta _{kj}}{(j-1)!}G^{(j-1)}(\lambda
_{k}),\qquad G(\lambda )\in
\mathbb{C}
_{2N}[\lambda ].  \label{5.12}
\end{equation}%
Then formula (\ref{5.10}) can be written in the form%
\begin{equation}
f_{n}=\left\langle \Omega ,F(\lambda )P_{n}(\lambda )\right\rangle ,\qquad n\in \{0,1,\ldots ,N-1\}.  \label{5.13}
\end{equation}%
From here by (\ref{5.11}) and the arbitrariness of $\{f_{m}\}_{m=0}^{N-1}$
it follows that the ``orthogonality''
relation%
\begin{equation}
\left\langle \Omega ,P_{m}(\lambda )P_{n}(\lambda )\right\rangle =\delta
_{mn},\qquad m,n\in \{0,1,\ldots ,N-1\},  \label{5.14}
\end{equation}%
holds. Further, in virtue of (\ref{5.8}) and (\ref{5.12}) we have also%
\begin{equation}
\left\langle \Omega ,P_{m}(\lambda )P_{N}(\lambda )\right\rangle =0,\qquad m\in \{0,1,\ldots ,N\}.  \label{5.15}
\end{equation}%
These mean by Theorem \ref{Th2.2} that the generalized spectral function of
the matrix $J$ has the form~(\ref{5.12}).
\end{proof}

\begin{definition}
The collection of the quantities%
\begin{equation*}
\{\lambda _{k},\beta _{kj}\ (j=1,\ldots ,m_{k},k=1,\ldots ,p)\},
\end{equation*}%
determining the structure of the generalized spectral function of the matrix
$J$ according to Theorem \ref{Th5.1}, we call the spectral data of the
matrix $J.$ For each $k\in \{1,\ldots ,p\}$ the sequence%
\begin{equation*}
\{\beta _{k1},\ldots ,\beta _{km_{k}}\}
\end{equation*}%
we call the normalizing chain (of the matrix $J$) associated with the
eigenvalue $\lambda _{k}$ (the sense of ``normalizing'' will be clear below in Section~\ref{section8}).
\end{definition}

If we delete the f\/irst row and the f\/irst column of the matrix $J$ given in (\ref{1.1}), then we get the new matrix%
\begin{equation*}
J^{(1)}=\left[
\begin{array}{ccccccc}
b_{0}^{(1)} & a_{0}^{(1)} & 0 & \cdots & 0 & 0 & 0 \\
a_{0}^{(1)} & b_{1}^{(1)} & a_{1}^{(1)} & \cdots & 0 & 0 & 0 \\
0 & a_{1}^{(1)} & b_{2}^{(1)} & \cdots & 0 & 0 & 0 \\
\vdots & \vdots & \vdots & \ddots & \vdots & \vdots & \vdots \\
0 & 0 & 0 & \ldots & b_{N-4}^{(1)} & a_{N-4}^{(1)} & 0 \\
0 & 0 & 0 & \cdots & a_{N-4}^{(1)} & b_{N-3}^{(1)} & a_{N-3}^{(1)} \\
0 & 0 & 0 & \cdots & 0 & a_{N-3}^{(1)} & b_{N-2}^{(1)}%
\end{array}%
\right] ,
\end{equation*}%
where%
\begin{gather*}
a_{n}^{(1)}=a_{n+1},\qquad n\in \{0,1,\ldots ,N-3\},
\\
b_{n}^{(1)}=b_{n+1},\qquad n\in \{0,1,\ldots ,N-2\}.
\end{gather*}%
The matrix $J^{(1)}$ is called the \textit{first truncated matrix }(with
respect to the matrix $J$).

\begin{theorem}
\label{Th5.2} The normalizing numbers $\beta _{kj}$ of the matrix $J$ can be
calculated by decomposing the rational function%
\begin{equation*}
-\frac{\det (J^{(1)}-\lambda I)}{\det (J-\lambda I)}
\end{equation*}%
into partial fractions.
\end{theorem}

\begin{proof}
Let us denote the polynomials of the f\/irst and the second kinds,
corresponding to the matrix $J^{(1)},$ by $P_{n}^{(1)}(\lambda )$ and $%
Q_{n}^{(1)}(\lambda ),$ respectively. It is easily seen that%
\begin{gather}
P_{n}^{(1)}(\lambda )=a_{0}Q_{n+1}(\lambda ),\qquad n\in \{0,1,\ldots
,N-1\},  \label{5.16}
\\
Q_{n}^{(1)}(\lambda )=\frac{1}{a_{0}}\{(\lambda -b_{0})Q_{n+1}(\lambda
)-P_{n+1}(\lambda )\},\qquad n\in \{0,1,\ldots ,N-1\}.  \label{5.17}
\end{gather}%
Indeed, both sides of each of these equalities are solutions of the same
dif\/ference equation%
\begin{equation*}
a_{n-1}^{(1)}y_{n-1}+b_{n}^{(1)}y_{n}+a_{n}^{(1)}y_{n+1}=\lambda y_{n},\qquad n\in \{0,1,\ldots ,N-2\},\qquad a_{N-2}^{(1)}=1,
\end{equation*}%
and the sides coincide for $n=-1$ and $n=0.$ Therefore the equality holds by
the uniqueness theorem for solutions.

Consequently, taking into account Lemma \ref{Lem2.1} and using (\ref{5.16}),
we have%
\begin{gather*}
\det (J^{(1)}-\lambda I)  = (-!)^{N-1}a_{0}^{(1)}a_{1}^{(1)}\cdots
a_{N-3}^{(1)}P_{N-1}^{(1)}(\lambda )
 = (-1)^{N-1}a_{1}\cdots a_{N-2}a_{0}Q_{N}(\lambda ).
\end{gather*}%
Comparing this with (\ref{2.6}), we get%
\begin{equation*}
\frac{Q_{N}(\lambda )}{P_{N}(\lambda )}=-\frac{\det (J^{(1)}-\lambda I)}{%
\det (J-\lambda I)}
\end{equation*}%
so that the statement of the theorem follows from (\ref{5.9}).
\end{proof}

\section{Inverse problem from the spectral data}\label{section6}

By the inverse spectral problem is meant the problem of recovering matrix $J$%
, i.e. its entries $a_{n}$ and $b_{n},$ from the spectral data.

\begin{theorem}
\label{Th6.1} Let  an arbitrary collection of complex numbers%
\begin{equation}
\{\lambda _{k},\beta _{kj}\ (j=1,\ldots ,m_{k},k=1,\ldots ,p)\}
\label{6.1}
\end{equation}%
be given, where $\lambda _{1},\lambda _{2},\ldots ,\lambda _{p}$ $(1\leq
p\leq N)$ are distinct, $1\leq m_{k}\leq N$, and $m_{1}+\dots +m_{p}=N$. In
order for this collection to be the spectral data for some Jacobi matrix $J$
of the form~\eqref{1.1} with entries belonging to the class~\eqref{1.2}, it
is necessary and sufficient that the following two conditions be satisfied:
\begin{enumerate}\itemsep=0pt
\item[$(i)$] $\sum\limits_{k=1}^{p}\beta _{k1}=1$;

\item[$(ii)$] $D_{n}\neq 0$,   for  $n\in \{1,2,\ldots ,N-1\}$, and $D_{N}=0$, where $D_{n}$ is defined by \eqref{3.24} in which%
\begin{equation}
s_{l}=\sum_{k=1}^{p}\sum_{j=1}^{n_{kl}}\binom{l}{j-1}\beta _{kj}\lambda
_{k}^{l-j+1},  \label{6.2}
\end{equation}%
$n_{kl}=\min \{m_{k},l+1\}$, $\binom{l}{j-1}$ is a binomial coefficient.
\end{enumerate}
\end{theorem}

\begin{proof}
The necessity of conditions of the theorem follows from Theorem \ref{Th3.2}
because the generali\-zed spectral function of the matrix $J$ is def\/ined by
the spectral data according to formula (\ref{5.1}) and therefore the
quantity (\ref{6.2}) coincides with $\left\langle \Omega ,\lambda
^{l}\right\rangle .$ Besides,%
\begin{equation*}
\sum_{k=1}^{p}\beta _{k1}=\left\langle \Omega ,1\right\rangle =s_{0}=D_{0}.
\end{equation*}%
Note that the condition $(iii)$ of Theorem \ref{Th3.1} holds with%
\begin{equation}
T(\lambda )=(\lambda -\lambda _{1})^{m_{1}}\cdots (\lambda -\lambda
_{p})^{m_{p}}.  \label{6.3}
\end{equation}

Let us prove the suf\/f\/iciency. Assume that we have a collection of quantities
(\ref{6.1}) satisfying the conditions of the theorem. Using these data we
construct the functional $\Omega $ on $%
\mathbb{C}
_{2N}[\lambda ]$ by formula~(\ref{5.1}). Then this functional $\Omega $
satisf\/ies the conditions of Theorem \ref{Th3.2} and therefore there exists a
matrix $J$ of the form (\ref{1.1}), (\ref{1.2}) for which $\Omega $ is the
generalized spectral function. Now we have to prove that the collection (\ref%
{6.1}) is the spectral data for the recovered matrix $J.$ For this purpose
we def\/ine the polynomials $P_{-1}(\lambda ),P_{0}(\lambda ),\ldots
,P_{N}(\lambda )$ as the solution of equation~(\ref{5.2}), constructed by
means of the matrix $J,$ under the initial conditions~(\ref{5.3}). Then the
rela\-tions~(\ref{2.7}),~(\ref{2.8}) and the equalities%
\begin{gather}
a_{n}=\left\langle \Omega ,\lambda P_{n}(\lambda )P_{n+1}(\lambda
)\right\rangle ,\qquad n\in \{0,1,\ldots ,N-2\},  \label{6.4}
\\
b_{n}=\left\langle \Omega ,\lambda P_{n}^{2}(\lambda )\right\rangle ,\qquad n\in \{0,1,\ldots ,N-1\}  \label{6.5}
\end{gather}%
hold. We show that (\ref{5.8}) holds which will mean, in particular, that $%
\lambda _{1},\ldots ,\lambda _{p}$ are eigenvalues of the matrix $J$ with
the multiplicities $m_{1},\ldots ,m_{p},$ respectively.

Let $T(\lambda )$ be def\/ined by (\ref{6.3}). Let us show that there exists a
constant $c$ such that%
\begin{equation}
a_{N-2}P_{N-2}(\lambda )+b_{N-1}P_{N-1}(\lambda )+cT(\lambda )=\lambda
P_{N-1}(\lambda )  \label{6.6}
\end{equation}%
for all $\lambda \in
\mathbb{C}
.$ If we prove this, then from here and (\ref{5.2}) with $%
y_{k}=P_{k}(\lambda )$ and $n=N-1$ we get that $P_{N}(\lambda )=cT(\lambda
). $

Since $\deg P_{n}(\lambda )=n$ $(0\leq n\leq N-1),$ $\deg T(\lambda
)=m_{1}+\dots +m_{p}=N,$ the polynomials $P_{0}(\lambda ),\ldots
,P_{N-1}(\lambda ),T(\lambda )$ form a basis of the linear space of all
polynomials of degree $\leq N.$ Therefore we have the decomposition%
\begin{equation}
\lambda P_{N-1}(\lambda )=cT(\lambda )+\sum_{n=0}^{N-1}c_{n}P_{n}(\lambda ),
\label{6.7}
\end{equation}%
where $c,c_{0},c_{1},\ldots ,c_{N-1}$ are some constants. By (\ref{6.3}) and
(\ref{5.1}) it follows that%
\begin{equation*}
\left\langle \Omega ,T(\lambda )P_{n}(\lambda )\right\rangle =0,\qquad
n\in \{0,1,\ldots ,N\}.
\end{equation*}%
Hence taking into account the relations (\ref{2.7}), (\ref{2.8}) and (\ref%
{6.4}), (\ref{6.5}), we f\/ind from (\ref{6.7}) that%
\begin{equation*}
c_{n}=0\quad (0\leq n\leq N-3),\qquad c_{N-2}=a_{N-2},\qquad
c_{N-1}=b_{N-1}.
\end{equation*}%
So (\ref{6.6}) is shown.

It remains to show that for each $\ k\in \{1,\ldots ,p\}$ the sequence $%
\{\beta _{k1},\ldots ,\beta _{km_{k}}\}$ is the normali\-zing chain of the
matrix $J$ associated with the eigenvalue $\lambda _{k}.$ Since we have
already shown that $\lambda _{k}$ is an eigenvalue of the matrix $J$ of the
multiplicity $m_{k},$ the normalizing chain of $J$ associated with the
eigenvalue $\lambda _{k}$ has the form $\{\widetilde{\beta }_{k1},\ldots ,%
\widetilde{\beta }_{km_{k}}\}.$ Therefore for $\left\langle \Omega
,G(\lambda )\right\rangle $ we have an equality of the form (\ref{5.1}) in
which $\beta _{kj}$ is replaced by $\widetilde{\beta }_{kj}.$ Subtracting
these two equalities for $\left\langle \Omega ,G(\lambda )\right\rangle $
each from other we get that%
\begin{equation*}
\sum_{k=1}^{p}\sum_{j=1}^{m_{k}}\frac{\beta _{kj}-\widetilde{\beta }_{kj}}{%
(j-1)!}G^{(j-1)}(\lambda _{k})=0\qquad \text{for all}\quad G(\lambda )\in
\mathbb{C}
_{2N}[\lambda ].
\end{equation*}%
Since the values $G^{(j-1)}(\lambda _{k})$ can be arbitrary numbers, we get
that $\beta _{kj}=\widetilde{\beta }_{kj}$ for all $k$ and $j.$
\end{proof}

Under the conditions of Theorem \ref{Th6.1} the entries $a_{n}$ and $b_{n}$
of the matrix $J$ for which the collection (\ref{6.1}) is spectral data, are
recovered by formulas (\ref{3.35}), (\ref{3.36}).

\section[Characterization of generalized spectral functions of real Jacobi matrices]{Characterization of generalized spectral functions\\ of real Jacobi matrices}\label{section7}

In this section, we characterize generalized spectral functions of real
Jacobi matrices among the generalized spectral functions of complex Jacobi
matrices. Let $m$ be a nonnegative integer. Denote by $%
\mathbb{R}
_{2m}[\lambda ]$ the ring of all polynomials in $\lambda $ of degree $\leq
2m $ with real coef\/f\/icients.

\begin{definition}
A linear functional $\Omega $ def\/ined on the space $%
\mathbb{C}
_{2m}[\lambda ]$ is said to be positive if%
\begin{equation*}
\left\langle \Omega ,G(\lambda )\right\rangle >0
\end{equation*}%
for all polynomials $G(\lambda )\in
\mathbb{R}
_{2m}[\lambda ],$ which are not identically zero and which satisfy the
inequality%
\begin{equation*}
G(\lambda )\geq 0,\qquad -\infty <\lambda <\infty .
\end{equation*}
\end{definition}

\begin{lemma}
\label{Lem7.1} If $\Omega $ is a positive functional on $%
\mathbb{C}
_{2m}[\lambda ],$ then it takes only real values on $%
\mathbb{R}
_{2m}[\lambda ].$
\end{lemma}

\begin{proof}
Since the functional $\Omega $ is positive, the values $\langle \Omega
,\lambda ^{2k}\rangle ,$ $k\in \{0,1,\ldots ,m\}$ are real (moreover
they are positive). Next, the monomial $\lambda ^{2k-1},$ $k\in \{1,2,\ldots
,m\}$ is represented as a dif\/ference of two nonnegative polynomials of
degree $2k:$%
\begin{equation*}
2\lambda ^{2k-1}=\lambda ^{2k-2}(\lambda +1)^{2}-\lambda ^{2k-2}(\lambda
^{2}+1).
\end{equation*}%
Therefore the values $ \langle \Omega ,\lambda ^{2k-1} \rangle ,$ $%
k\in \{1,2,\ldots ,m\}$ are also real to be a dif\/ference of two positive
numbers. Thus, $\left\langle \Omega ,\lambda ^{n}\right\rangle $ is real for
any $n\in \{0,1,\ldots ,2m\}.$ Hence $\left\langle \Omega ,G(\lambda
)\right\rangle $ is real for any $G(\lambda )\in
\mathbb{R}
_{2m}[\lambda ].$
\end{proof}

\begin{lemma}
\label{Lem7.2} A linear functional $\Omega $ on $%
\mathbb{C}
_{2m}[\lambda ]$ is positive if and only if $D_{n}>0$ for all $n\in
\{0,1,\ldots ,m\}$, where%
\begin{equation*}
D_{n}=\left\vert
\begin{array}{cccc}
s_{0} & s_{1} & \cdots & s_{n} \\
s_{1} & s_{2} & \cdots & s_{n+1} \\
\vdots & \vdots & \ddots & \vdots \\
s_{n} & s_{n+1} & \cdots & s_{2n}%
\end{array}%
\right\vert ,\qquad n=0,1,\ldots ,m,
\end{equation*}%
in which%
\begin{equation*}
s_{l}= \langle \Omega ,\lambda ^{l} \rangle ,\qquad
l=0,1,\ldots ,2m.
\end{equation*}
\end{lemma}

\begin{proof}
Any polynomial $G(\lambda )\in
\mathbb{R}
_{2m}[\lambda ]$ which is not identically zero and which satisf\/ies the
inequality%
\begin{equation}
G(\lambda )\geq 0,\qquad -\infty <\lambda <\infty ,  \label{7.1}
\end{equation}%
can be represented in the form%
\begin{equation}
G(\lambda )=[A(\lambda )]^{2}+[B(\lambda )]^{2},  \label{7.2}
\end{equation}%
where $A(\lambda )$, $B(\lambda )$ are polynomials of degrees $\leq m$ with
real coef\/f\/icients. Indeed, it follows from~(\ref{7.1}) that the polynomial $%
G(\lambda )$ has even degree: $\deg G(\lambda )=2p,$ where $p\leq m.$
Therefore its decomposition into linear factors has the form%
\begin{equation*}
G(\lambda )=c\prod_{k=1}^{p}(\lambda -\alpha _{k}-i\beta
_{k})(\lambda -\alpha _{k}+i\beta _{k}),
\end{equation*}%
where $c>0,$ $\beta _{k}\geq 0,$ $\alpha _{k}$ are real (among the roots $%
\alpha _{k}+i\beta _{k},$ of course, may be equal). Now setting%
\begin{equation*}
\sqrt{c}\prod_{k=1}^{p}(\lambda -\alpha _{k}-i\beta
_{k})=A(\lambda )+iB(\lambda ),
\end{equation*}%
we get   that the polynomials $A(\lambda ),$ $B(\lambda )$ have real
coef\/f\/icients and (\ref{7.2}) holds.

Now writing%
\begin{equation*}
A(\lambda )=\sum_{k=1}^{p}x_{k}\lambda ^{k},\qquad B(\lambda
)=\sum_{k=1}^{p}y_{k}\lambda ^{k},
\end{equation*}%
where $x_{k}$, $y_{k}$ are real numbers, we f\/ind%
\begin{equation*}
\left\langle \Omega ,G(\lambda )\right\rangle
=\sum_{j,k=0}^{p}s_{j+k}x_{j}x_{k}+\sum_{j,k=0}^{p}s_{j+k}y_{j}y_{k}.
\end{equation*}%
This implies the statement of the lemma.
\end{proof}

\begin{theorem}
\label{Th7.3} In order for a given linear functional $\Omega $ on $%
\mathbb{C}
_{2N}[\lambda ]$, to be the generalized spectral function for a real Jacobi
matrix of the form \eqref{1.1}, \eqref{1.3} it is necessary and sufficient
that the following three conditions be satisfied:
\begin{enumerate}\itemsep=0pt
\item[$(i)$] $\left\langle \Omega ,1\right\rangle =1$;

\item[$(ii)$] $\Omega $ is positive on $%
\mathbb{C}
_{2N-2}[\lambda ]$;

\item[$(iii)$] there exists a polynomial $T(\lambda )$ of degree $N$ such that
$\left\langle \Omega ,G(\lambda )T(\lambda )\right\rangle =0$ for all
polynomials $G(\lambda )$ with $\deg G(\lambda )\leq N.$
\end{enumerate}
\end{theorem}

\begin{proof}
\textit{Necessity}. The condition $\left\langle \Omega ,1\right\rangle =1$
follows from (\ref{2.7}) with $m=n=0.$ To prove positivity on $%
\mathbb{C}
_{2N-2}[\lambda ]$ of the generalized spectral function $\Omega $ of the
real Jacobi matrix $J,$ take an arbitrary polynomial $G(\lambda )\in
\mathbb{R}
_{2N-2}[\lambda ]$ which is not identically zero and which satisf\/ies the
inequality%
\begin{equation*}
G(\lambda )\geq 0,\qquad -\infty <\lambda <\infty .
\end{equation*}%
This polynomial can be represented in the form (see the proof of Lemma \ref%
{Lem7.2})%
\begin{equation}
G(\lambda )=[A(\lambda )]^{2}+[B(\lambda )]^{2},  \label{7.3}
\end{equation}%
where $A(\lambda ),$ $B(\lambda )$ are polynomials of degrees $\leq N-1$
with real coef\/f\/icients. Since the polynomials $P_{0}(\lambda ),P_{1}(\lambda
),\ldots ,P_{N-1}(\lambda )$ have real coef\/f\/icients (because $J$ is a real
matrix) and they form a basis of $%
\mathbb{R}
_{N-1}[\lambda ],$ we can write the decompositions%
\begin{equation*}
A(\lambda )=\sum_{k=1}^{N-1}c_{k}P_{k}(\lambda ),\qquad B(\lambda
)=\sum_{k=1}^{N-1}d_{k}P_{k}(\lambda ),
\end{equation*}%
where $c_{k}$, $d_{k}$ are real numbers not all zero. Therefore using the
``orthogonality'' property (\ref{2.7}) we
get from (\ref{7.3}),%
\begin{equation*}
\left\langle \Omega ,G(\lambda )\right\rangle
=\sum_{j,k=0}^{N-1}(c_{k}^{2}+d_{k}^{2})>0.
\end{equation*}

The property of $\Omega $ indicated in the condition $(iii)$ of the theorem
follows from (\ref{2.8}) if we take $T(\lambda )=P_{N}(\lambda ).$

\textit{Sufficiency.} It follows from the conditions of the theorem that all
the conditions of Theorem~\ref{Th3.1} are satisf\/ied. In fact, we need to
verify only the condition $(ii)$ of Theorem \ref{Th3.1}. Let for some
polynomial $G(\lambda ),$ $\deg G(\lambda )=n\leq N-1,$%
\begin{equation}
\left\langle \Omega ,G(\lambda )H(\lambda )\right\rangle =0  \label{7.4}
\end{equation}%
for all polynomials $H(\lambda )$, $\deg H(\lambda )=n.$ We have to show
that then $G(\lambda )\equiv 0.$ Setting
\begin{equation*}
G(\lambda )=\sum_{k=0}^{n}g_{k}\lambda ^{k},\qquad H(\lambda
)=\sum_{j=0}^{n}h_{j}\lambda ^{j},
\end{equation*}%
we get from (\ref{7.4}) that%
\begin{equation*}
\sum_{j=0}^{n}h_{j}\left( \sum_{k=0}^{n}g_{k}s_{k+j}\right) =0.
\end{equation*}%
Since $h_{0},h_{1},\ldots ,h_{n}$ $(h_{n}\neq 0)$ are arbitrary, the last
equation gives%
\begin{equation}
\sum_{k=0}^{n}g_{k}s_{k+j}=0,\qquad j=0,1,\ldots ,n.  \label{7.5}
\end{equation}%
This is a linear homogeneous system of algebraic equations with respect to $%
g_{0},g_{1},\ldots ,g_{n}$ and the determinant of this system coincides with
the determinant $D_{n}.$ From the condition $(ii)$ of the theorem it follows
by Lemma~\ref{Lem7.2} that $D_{n}>0.$ So $D_{n}\neq 0$ and hence system (\ref{7.5}) has only the trivial solution $g_{0}=g_{1}=\dots =g_{n}.$

Thus, all the conditions of Theorem~\ref{Th3.1} are satisf\/ied. Therefore
there exists, generally speaking, a complex Jacobi matrix $J$ of the form (\ref{1.1}), (\ref{1.2}) for which $\Omega $ is the generalized spectral
function. This matrix $J$ is constructed by using formulas (\ref{3.35}), (\ref{3.36}). It remains to show that the matrix~$J$ is real. But this
follows from the fact that by Lemma~\ref{Lem7.1} and Lemma~\ref{Lem7.2} we
have $D_{n}>0$ for $n\in \{0,1,\ldots ,N-1\}$ and the determinants $\Delta
_{n}$ are real. Therefore formulas (\ref{3.35}), (\ref{3.36}) imply that the
matrix $J$ is real.
\end{proof}

If we take into account Lemma~\ref{Lem7.2}, then it is easily seen from the
proof of Theorem~\ref{Th3.2} that Theorem~\ref{Th7.3} is equivalent to the
following theorem.

\begin{theorem}
\label{Th7.4} In order for a given linear functional $\Omega $, defined on $%
\mathbb{C}
_{2N}[\lambda ]$, to be the generalized spectral function for some real
Jacobi matrix $J$ of the form \eqref{1.1} with entries belonging to the
class \eqref{1.3}, it is necessary and sufficient that
\begin{equation*}
D_{0}=1,\qquad D_{n}>0\quad (n=1,2,\ldots ,N-1),\qquad \text{and} \qquad
D_{N}=0,
\end{equation*}%
where $D_{n}$ is defined by \eqref{3.24} and \eqref{3.23}.
\end{theorem}

Under the conditions of Theorem \ref{Th7.4} the entries $a_{n}$ and $b_{n}$
of the matrix $J$ for which the functional $\Omega $ is the spectral
function are recovered by formulas (\ref{3.35}), (\ref{3.36}).

\section[Structure of generalized spectral functions of real Jacobi matrices]{Structure of generalized spectral functions\\ of real Jacobi matrices}\label{section8}

First we prove two lemmas which hold for any complex Jacobi matrix $J$ of
the form~(\ref{1.1}),~(\ref{1.2}). Having the matrix $J$ consider the
dif\/ference equation (\ref{5.2}) and let $\{P_{n}(\lambda )\}_{n{=}-1}^{N} $ and
$\{Q_{n}(\lambda )\}_{n{=}-1}^{N}\!$ be the solutions of this equation
satisfying the initial conditions (\ref{5.3}) and (\ref{5.4}), respectively.

\begin{lemma}
\label{Lem8.1} The equation%
\begin{equation}
P_{N-1}(\lambda )Q_{N}(\lambda )-P_{N}(\lambda )Q_{N-1}(\lambda )=1
\label{8.1}
\end{equation}%
holds.
\end{lemma}

\begin{proof}
Multiply f\/irst of the equations%
\begin{gather}
a_{n-1}P_{n-1}(\lambda )+b_{n}P_{n}(\lambda )+a_{n}P_{n+1}(\lambda )=\lambda
P_{n}(\lambda ),  \label{8.2}
\\
n\in \{0,1,\ldots ,N-1\},\qquad a_{-1}=a_{N-1}=1,\nonumber
\\
a_{n-1}Q_{n-1}(\lambda )+b_{n}Q_{n}(\lambda )+a_{n}Q_{n+1}(\lambda )=\lambda
Q_{n}(\lambda ),\nonumber
\\
n\in \{0,1,\ldots ,N-1\},\qquad a_{-1}=a_{N-1}=1,\nonumber
\end{gather}%
by $Q_{n}(\lambda )$ and second by $P_{n}(\lambda )$ and subtract the second
result from the f\/irst to get%
\begin{gather*}
a_{n-1}[P_{n-1}(\lambda )Q_{n}(\lambda )-P_{n}(\lambda )Q_{n-1}(\lambda )]
\\
\qquad{}=a_{n}[P_{n}(\lambda )Q_{n+1}(\lambda )-P_{n+1}(\lambda )Q_{n}(\lambda )],%
\qquad n\in \{0,1,\ldots ,N-1\}.
\end{gather*}%
This means that the expression (Wronskian of the solutions $P_{n}(\lambda )$
and $Q_{n}(\lambda )$)%
\begin{equation*}
a_{n}[P_{n}(\lambda )Q_{n+1}(\lambda )-P_{n+1}(\lambda )Q_{n}(\lambda )]
\end{equation*}%
does not depend on $n\in \{-1,0,1,\ldots ,N-1\}.$ On the other hand the
value of this expression at $n=-1$ is equal to $1$ by (\ref{5.3}), (\ref{5.4}%
), and $a_{-1}=1.$ Therefore%
\begin{equation*}
a_{n}[P_{n}(\lambda )Q_{n+1}(\lambda )-P_{n+1}(\lambda )Q_{n}(\lambda )]=1\qquad
\text{for all} \quad n\in \{-1,0,1,\ldots ,N-1\}.
\end{equation*}%
Putting here, in particular, $n=N-1,$ we arrive at (\ref{8.1}).
\end{proof}

\begin{lemma}
\label{Lem8.2} The equation%
\begin{equation}
P_{N-1}(\lambda )P_{N}^{\prime }(\lambda )-P_{N}(\lambda )P_{N-1}^{\prime
}(\lambda )=\sum_{n=0}^{N-1}P_{n}^{2}(\lambda )  \label{8.3}
\end{equation}%
holds, where the prime denotes the derivative with respect to $\lambda $.
\end{lemma}

\begin{proof}
Dif\/ferentiating equation (\ref{8.2}) with respect to $\lambda ,$ we get%
\begin{gather}
a_{n-1}P_{n-1}^{\prime }(\lambda )+b_{n}P_{n}^{\prime }(\lambda
)+a_{n}P_{n+1}^{\prime }(\lambda )=\lambda P_{n}^{\prime }(\lambda
)+P_{n}(\lambda ),  \label{8.4}
\\
n\in \{0,1,\ldots ,N-1\},\qquad a_{-1}=a_{N-1}=1.\nonumber
\end{gather}%
Multiplying equation (\ref{8.2}) by $P_{n}^{\prime }(\lambda )$ and equation
(\ref{8.4}) by $P_{n}(\lambda ),$ and subtracting the left and right members
of the resulting equations, we get
\begin{gather*}
a_{n-1}[P_{n-1}(\lambda )P_{n}^{\prime }(\lambda )-P_{n}(\lambda
)P_{n-1}^{\prime }(\lambda )]-a_{n}[P_{n}(\lambda )P_{n+1}^{\prime }(\lambda
)-P_{n+1}(\lambda )P_{n}^{\prime }(\lambda )]=-P_{n}^{2}(\lambda ),
\\
\qquad {} n\in \{0,1,\ldots ,N-1\}.
\end{gather*}%
Summing the last equation for the values $n=0,1,\ldots ,m$ ($m\leq N-1$) and
using the initial conditions (\ref{5.3}), we obtain%
\begin{equation*}
a_{m}[P_{m}(\lambda )P_{m+1}^{\prime }(\lambda )-P_{m+1}(\lambda
)P_{m}^{\prime }(\lambda )]=\sum_{n=0}^{m}P_{n}^{2}(\lambda ),\qquad
m\in \{0,1,\ldots ,N-1\}.
\end{equation*}%
Setting here, in particular, $m=N-1$ and taking into account $a_{N-1}=1,$ we
get (\ref{8.3}).
\end{proof}

Now we consider real Jacobi matrices of the form (\ref{1.1}), (\ref{1.3}).

\begin{lemma}
\label{Lem8.3} For any real Jacobi matrix $J$ of the form \eqref{1.1}, \eqref{1.3} the roots of the polyno\-mial~$P_{N}(\lambda )$ are simple.
\end{lemma}

\begin{proof}
Let $\lambda _{0}$ be a root of the polynomial $P_{N}(\lambda ).$ The root $%
\lambda _{0}$ is an eigenvalue of the matrix $J$ by (\ref{2.6}) and hence it
is real by Hermiticity of $J.$ Putting $\lambda =\lambda _{0}$ in (\ref{8.3}%
) and using $P_{N}(\lambda _{0})=0,$ we get%
\begin{equation}
P_{N-1}(\lambda _{0})P_{N}^{\prime }(\lambda
_{0})=\sum_{n=0}^{N-1}P_{n}^{2}(\lambda _{0}).  \label{8.5}
\end{equation}%
The right-hand side of (\ref{8.5}) is dif\/ferent from zero because the
polynomials $P_{n}(\lambda )$ have real coef\/f\/icients and hence are real for
real values of $\lambda ,$ and $P_{0}(\lambda )=1.$ Consequently $%
P_{N}^{\prime }(\lambda _{0})\neq 0,$ that is, the root $\lambda _{0}$ of
the polynomial $P_{N}(\lambda )$ is simple.
\end{proof}

\begin{lemma}
\label{Lem8.4} Any real Jacobi matrix $J$ of the form \eqref{1.1}, \eqref{1.3} has precisely $N$ real and distinct eigenvalues.
\end{lemma}

\begin{proof}
The reality of eigenvalues of $J$ follows from its Hermiticity. Next, the
eigenvalues of $J$ coincide, by (\ref{2.6}), with the roots of the
polynomial $P_{N}(\lambda ).$ This polynomial of degree $N$ has $N$ roots.
These roots are pairwise distinct by Lemma \ref{Lem8.3}.
\end{proof}

The following theorem describes the structure of generalized spectral
functions of real Jacobi matrices.

\begin{theorem}
\label{Th8.4} Let $J$ be a real Jacobi matrix of the form \eqref{1.1}, \eqref{1.3} and $\Omega $ be its generalized spectral function. Then for any $%
G(\lambda )\in
\mathbb{C}
_{2N}[\lambda ]$%
\begin{equation}
\left\langle \Omega ,G(\lambda )\right\rangle =\sum_{k=1}^{N}\beta
_{k}G(\lambda _{k}),  \label{8.6}
\end{equation}%
where $\lambda _{1},\ldots ,\lambda _{N}$ are the eigenvalues of the matrix~$J$ and $\beta _{1},\ldots ,\beta _{N}$ are positive real numbers uniquely
determined by the matrix~$J$.
\end{theorem}

\begin{proof}
By Lemma \ref{Lem8.3}, the roots $\lambda _{1},\ldots ,\lambda _{N}$ of the
polynomial $P_{N}(\lambda )$ are simple. Therefore the formula (\ref{8.6})
follows from (\ref{5.1}) and the decomposition (\ref{5.9}) takes the form%
\begin{equation*}
\frac{Q_{N}(\lambda )}{P_{N}(\lambda )}=\sum_{k=1}^{N}\frac{\beta _{k}}{%
\lambda -\lambda _{k}}.
\end{equation*}%
Hence%
\begin{equation}
Q_{N}(\lambda _{k})=\beta _{k}P_{N}^{\prime }(\lambda _{k}).  \label{8.7}
\end{equation}%
On the other hand, putting $\lambda =\lambda _{k}$ in (\ref{8.1}) and (\ref%
{8.3}) and taking into account that $P_{N}(\lambda _{k})=0,$ we get%
\begin{gather}
P_{N-1}(\lambda _{k})Q_{N}(\lambda _{k})=1,  \label{8.8}
\\
P_{N-1}(\lambda _{k})P_{N}^{\prime }(\lambda
_{k})=\sum_{n=0}^{N-1}P_{n}^{2}(\lambda _{k}),  \label{8.9}
\end{gather}%
respectively. Comparing (\ref{8.7}), (\ref{8.8}), and (\ref{8.9}), we f\/ind
that%
\begin{equation}
\beta _{k}=\left\{ \sum_{n=0}^{N-1}P_{n}^{2}(\lambda _{k})\right\} ^{-1},
\label{8.10}
\end{equation}%
whence we get, in particular, that $\beta _{k}>0.$
\end{proof}

Since $\{P_{n}(\lambda _{k})\}_{n=0}^{N-1}$ is an eigenvector of the matrix $%
J$ corresponding to the eigenvalue $\lambda _{k},$ it is natural, according
to the formula (\ref{8.10}), to call the numbers $\beta _{k}$ the \textit{%
normalizing numbers} of the matrix $J.$

\begin{definition}
The collection of the eigenvalues and normalizing numbers%
\begin{equation*}
\{\lambda _{k},\beta _{k}\ (k=1,\ldots ,N)\}
\end{equation*}%
of the matrix $J$ of the form (\ref{1.1}), (\ref{1.3}) we call the spectral
data of this matrix.
\end{definition}

\begin{remark}
Assuming that $\lambda _{1}<\lambda _{2}<\dots <\lambda _{N},$ let us
introduce the nondecreasing step function $\omega (\lambda )$ on $(-\infty
,\infty )$ by%
\begin{equation*}
\omega (\lambda )=\sum_{\lambda _{k}\leq \lambda }\beta _{k},
\end{equation*}%
where $\omega (\lambda )=0$ if there is no\ $\lambda _{k}\leq \lambda .$ So
the eigenvalues of the matrix $J$ are the points of increase of the function
$\omega (\lambda ).$ Then equality (\ref{8.6}) can be written as%
\begin{equation*}
\left\langle \Omega ,G(\lambda )\right\rangle =\int_{-\infty }^{\infty
}G(\lambda )d\omega (\lambda ),
\end{equation*}%
where the integral is a Stieltjes integral. Therefore the orthogonality
relation (\ref{5.14}) can be written as%
\begin{equation*}
\int_{-\infty }^{\infty }P_{m}(\lambda )P_{n}(\lambda )d\omega (\lambda
)=\delta _{mn},\qquad m,n\in \{0,1,\ldots ,N-1\}
\end{equation*}%
and the expansion formula (\ref{5.13}) as%
\begin{equation*}
f_{n}=\int_{-\infty }^{\infty }F(\lambda )P_{n}(\lambda )d\omega (\lambda ),%
\qquad n\in \{0,1,\ldots ,N-1\},
\end{equation*}%
where $F(\lambda )$ is def\/ined by (\ref{5.11}):%
\begin{equation*}
F(\lambda )=\sum_{m=0}^{N-1}f_{m}P_{m}(\lambda ).
\end{equation*}%
Such function $\omega (\lambda )$ is known as a spectral function (see,
e.g., \cite{[21]}) of the operator (matrix) $J.$ This explains the source of
the term ``generalized spectral function'' used in the complex case.
\end{remark}

\section{Inverse spectral problem for real Jacobi matrices}\label{section9}

By the inverse spectral problem  for real Jacobi matrices we mean the
problem of recovering the matrix, i.e.\ its entries, from the spectral data.

\begin{theorem}
\label{Th9.1} Let   an arbitrary collection of numbers%
\begin{equation}
\{\lambda _{k},\beta _{k}\ (k=1,\ldots ,N)\}  \label{9.1}
\end{equation}%
be given. In order for this collection to be the spectral data for a real
Jacobi matrix $J$ of the form~\eqref{1.1} with entries belonging to the
class \eqref{1.3}, it is necessary and sufficient that the following two
conditions be satisfied:
\begin{enumerate}\itemsep=0pt
\item[$(i)$] The numbers $\lambda _{1},\ldots ,\lambda _{N}$ are real and
distinct.

\item[$(ii)$] The numbers $\beta _{1},\ldots ,\beta _{N}$ are positive and
such that $\sum\limits_{k=1}^{N}\beta _{k}=1.$
\end{enumerate}
\end{theorem}

\begin{proof}
The necessity of the conditions of the theorem was proved above. To prove
the suf\/f\/i\-cien\-cy, assume that we have a collection of quantities (\ref{9.1})
satisfying the conditions of the theorem. Using these data we construct the
functional $\Omega $ on $%
\mathbb{C}
_{2N}[\lambda ]$ by the formula%
\begin{equation}
\left\langle \Omega ,G(\lambda )\right\rangle =\sum_{k=1}^{N}\beta
_{k}G(\lambda _{k}),\qquad G(\lambda )\in
\mathbb{C}
_{2N}[\lambda ].  \label{9.2}
\end{equation}%
Then this functional $\Omega $ satisf\/ies the conditions of Theorem~\ref{Th7.3}. Indeed, we have%
\begin{equation*}
\left\langle \Omega ,1\right\rangle =\sum_{k=1}^{N}\beta _{k}=1,
\end{equation*}%
Next, let $G(\lambda )\in
\mathbb{R}
_{2N-2}[\lambda ]$ be an arbitrary polynomial which is not identically zero
and which satisf\/ies the inequality%
\begin{equation*}
G(\lambda )\geq 0,\qquad -\infty <\lambda <\infty .
\end{equation*}%
This polynomial can be represented in the form (see the proof of Lemma \ref%
{Lem7.2})%
\begin{equation*}
G(\lambda )=[A(\lambda )]^{2}+[B(\lambda )]^{2},
\end{equation*}%
where $A(\lambda ),$ $B(\lambda )$ are polynomials of degrees $\leq N-1$
with real coef\/f\/icients. Then%
\begin{equation}
\left\langle \Omega ,G(\lambda )\right\rangle =\sum_{k=1}^{N}\beta
_{k}G(\lambda _{k})=\sum_{k=1}^{N}\beta _{k}[A(\lambda
_{k})]^{2}+\sum_{k=1}^{N}\beta _{k}[B(\lambda _{k})]^{2}\geq 0.  \label{9.3}
\end{equation}%
We have to show that the equality sign in (\ref{9.3}) is impossible. If we
have the equality in (\ref{9.3}), then, since all the $\beta _{k}$ are
positive, we get
\begin{equation*}
A(\lambda _{1})=\dots =A(\lambda _{N})=0\qquad \text{and} \qquad B(\lambda
_{1})=\dots =B(\lambda _{N})=0.
\end{equation*}%
Hence $A(\lambda )\equiv 0$ and $B(\lambda )\equiv 0$ because $\lambda
_{1},\ldots ,\lambda _{N}$ are distinct and $\deg A(\lambda )\leq N-1,$ $%
\deg B(\lambda )\leq N-1.$ Therefore we get $G(\lambda )\equiv 0$ which is a
contradiction. Finally, if we take%
\begin{equation*}
T(\lambda )=(\lambda -\lambda _{1})\cdots (\lambda -\lambda _{N}),
\end{equation*}%
then the condition $(iii)$ of Theorem \ref{Th7.3} is also satisf\/ied: for any
polynomial $G(\lambda ),$%
\begin{equation*}
\left\langle \Omega ,G(\lambda )T(\lambda )\right\rangle
=\sum_{k=1}^{N}\beta _{k}G(\lambda _{k})T(\lambda _{k})=0.
\end{equation*}%
Thus, the functional $\Omega $ def\/ined by the formula (\ref{9.2}) satisf\/ies
all the conditions of Theorem \ref{Th7.3}. Therefore there exists a real
Jacobi matrix $J$ of the form (\ref{1.1}), (\ref{1.3}) for which $\Omega $
is the generalized spectral function. Further, from the proof of suf\/f\/iciency
of the conditions of Theorem \ref{Th6.1} it follows that the collection $%
\{\lambda _{k},\beta _{k} \ (k=1,\ldots ,N)\}$ is the spectral data for the
recovered matrix $J.$
\end{proof}

Note that under the conditions of Theorem \ref{Th9.1} the entries $a_{n}$
and $b_{n}$ of the matrix $J$ for which the collection (\ref{9.1}) is
spectral data, are recovered by formulas (\ref{3.35}), (\ref{3.36}).

\subsection*{Acknowledgements}

This work was supported by Grant 106T549 from the Scientif\/ic and Technological Research Council of Turkey (TUBITAK).

\pdfbookmark[1]{References}{ref}
\LastPageEnding

\end{document}